\input psfig.sty
\input amssym.def
\input amssym.tex

\magnification=1200

\catcode`\@=11

\hsize=125 mm   \vsize =175mm \hoffset=4mm    \voffset=10mm
\pretolerance=500 \tolerance=1000 \brokenpenalty=5000

%
%
%
%
%
%

\newif\ifpagetitre      \pagetitretrue
\newtoks\hautpagetitre  \hautpagetitre={\hfil}
\newtoks\baspagetitre   \baspagetitre={\hfil}

\newtoks\auteurcourant  \auteurcourant={\hfil}
\newtoks\titrecourant   \titrecourant={\hfil}

\newtoks\hautpagegauche \newtoks\hautpagedroite
\hautpagegauche={\hfil\the\auteurcourant\hfil}
\hautpagedroite={\hfil\the\titrecourant\hfil}

\newtoks\baspagegauche  \baspagegauche={\hfil\tenrm\folio\hfil}
\newtoks\baspagedroite  \baspagedroite={\hfil\tenrm\folio\hfil}

\headline={\ifpagetitre\the\hautpagetitre
\else\ifodd\pageno\the\hautpagedroite
\else\the\hautpagegauche\fi\fi}

\footline={\ifpagetitre\the\baspagetitre
\global\pagetitrefalse
\else\ifodd\pageno\the\baspagedroite
\else\the\baspagegauche\fi\fi}

\hautpagetitre={\hfill\tenbf \hfill} \hautpagetitre={\hfill}
\hautpagegauche={\tenrm\folio\hfill\the\auteurcourant}
\hautpagedroite={\tenrm\the\titrecourant\hfill\folio}
\baspagegauche={\hfil} \baspagedroite={\hfil}
\auteurcourant{Marmi, Moussa, Yoccoz} \titrecourant{Affine i.e.m.\
with a wandering interval }
\hfuzz=0.3pt
\font\tit=cmb10 scaled \magstep1

\def\H{\Bbb H}
\def\R{\Bbb R}

\def\Z{\Bbb Z}
\def\Q{\Bbb Q}

\def\1{\Bbb I}
\def\N{{\Bbb N}}

\def\P{{\Bbb P}}
\def\hH+{\hat{\H}^{+}}
\def\hHZ+{\widehat{\H^{+}/\Z}}

\def\PGL2Z{\hbox{PGL}\, (2,\Z)}
\def\GL2Z{\hbox{GL}\, (2,\Z)}
\def\SL2Z{\hbox{SL}\, (2,\Z)}

\def\cA{{\cal A}}

\def\cD{{\cal D}}

\def\cR{{\cal R}}

\def\proof{\noindent{\it Proof.\ }}
\def\qed{\hfill$\square$\par\smallbreak}
\def\Proc#1#2\par{\medbreak \noindent {\bf #1\enspace }{\sl #2}%
\par\ifdim \lastskip <\medskipamount \removelastskip \penalty 55\medskip \fi}%
\def\Def#1#2\par{\medbreak \noindent {\bf #1\enspace }{#2}%
\par\ifdim \lastskip <\medskipamount \removelastskip \penalty 55\medskip
\fi}
\def\qed{\hfill$\square$\par\smallbreak}
\def\hfl#1#2{\smash{\mathop{\hbox to 12mm{\rightarrowfill}}
\limits^{\scriptstyle#1}_{\scriptstyle#2}}}

\def\build#1_#2^#3{\mathrel{\mathop{\kern 0pt#1}\limits_{#2}^{#3}}}



\centerline{\tit Affine interval exchange maps with a wandering
interval}
\bigskip \centerline{S. Marmi\footnote{$^1$}{Scuola Normale
Superiore, Piazza dei Cavalieri 7, 56126 Pisa, Italy}, P.
Moussa\footnote{$^2$}{Institut de Physique Th\'eorique, CEA/Saclay,
91191 Gif-Sur-Yvette, France} and J.-C.
Yoccoz\footnote{$^3$}{Coll\`ege de France, 3, Rue d'Ulm, 75005
Paris }} \vskip 1.5 truecm \centerline{\bf  Abstract}
\smallskip
For almost all interval exchange maps $T_0$, with combinatorics of
genus $g\ge 2$, we construct affine interval exchange maps $T$
which are semi--conjugate to $T_0$ and have a wandering interval.
\bigskip
Mathematical Review Classification: Primary: 37C15 (Topological and differentiable equivalence,
conjugacy, invariants, moduli, classification);  Secondary:
37E05 (maps of the interval), 11J70 (Continued fractions and generalizations)
\bigskip
\centerline{CONTENTS}
\smallskip{
\item{0.} Introduction
\item{1.}  The continued fraction algorithm for interval exchange
maps
\item{1.1} Interval exchange maps
\item{1.2} The elementary step of the Rauzy--Veech
algorithm
\item{1.3} Rauzy diagrams
\item{1.4} The Rauzy--Veech and Zorich algorithms
\item{1.5} Dynamics of the continued fraction algorithms
\item{1.6} The continued fraction algorithm for generalized i.e.m.\
\item{2.} Deformations of affine interval exchange maps
\item{2.1} The set $\hbox{Aff}^{(1)}\, (\underline\gamma , w)$
\item{2.2} Affine motions
\item{3.} Wandering intervals for affine interval exchange maps
\item{3.1} The Zorich cocycle
\item{3.2} Statement of the result
\item{3.3} Reduction to a statement on Birkhoff sums
\item{3.4} Limit shapes for Birkhoff sums
\item{3.5} On the direction of $w$
\item{3.6} Consequences for limit shapes
\item{3.7} Proof of the Proposition in 3.3.1
\item{} References

 }

\vfill\eject

\vskip .5 truecm \noindent {\bf 0. Introduction}

\vskip 1. truecm
Quasiperiodic systems play a very important role in the theory of dynamical systems and in mathematical physics.

Irrational rotations of the circle are the prototype of quasiperiodic dynamics. The suspension of  these rotations produces linear flows on the two-dimensional torus. When analyzing the recurrence of rotations or the suspended flows, the modular group $\hbox{GL}\, (2,\Z)$ is of fundamental importance, providing the renormalization scheme associated to the continuous fraction of the rotation number.

Poincar\'e proved that any orientation-preserving homeomorphism of the circle with no periodic orbit is semi-conjugate to an irrational rotation. Later Denjoy constructed examples of $C^r$ diffeomorphisms with irrational rotation number and a wandering interval if $r<2$. He also proved that any $C^2$ diffeomorphism with no periodic orbit is conjugate to an irrational rotation. Actually, this result is also true for piecewise-affine homeomorphisms [He].

A natural generalization of the linear flows on the two-dimensional torus is obtained by considering linear flows on compact surfaces of higher genus, called translation surfaces. By a Poincar\'e section their dynamics can be reduced to (standard) interval exchange maps (i.e.m.\ ), which generalize rotations of the circle.

Let $\cal A$ be an alphabet with
$d\ge 2$ elements. A (standard) i.e.m.\ $T$ on an interval $I$ (of finite length) is determined by two partitions $(I_a^t )$, $(I_a^b )$, of $I$ with $I_a^t$, $I_a^b$ of the same length, the restriction of $T$ to $I_a^t$  being a translation with image $I_a^b$ . Thus $T$ is orientation-preserving and preserves Lebesgue measure. By relaxing the requirement on the lengths and only asking that the restriction of $T$ to $I_a^t$ is an orientation-preserving homeomorphism onto $I_a^b$ one obtains the definition of a generalized i.e.m. A special class of generalized i.e.m.\ , namely affine i.e.m.\ are considered in this paper: we require that the restriction of $T$ to $I_a^t$ is affine (and orientation-preserving). When $d=2$, by identifying the endpoints of $I$, standard i.e.m.\ correspond to rotations of the circle and generalized i.e.m.\ to homeomorphisms of the circle.

The ordering of the subintervals in the two partitions of $I$ constitute the combinatorial data for the i.e.m.\ $T$. One says that a standard i.e.m. has no connection if every orbit can be extended indefinitely in the future or in the past (or both) without going through the endpoints of the subintervals; Keane [Ke] has shown that such an i.e.m. is minimal. When $d=2$, this corresponds exactly to irrational rotations.

Following Rauzy [Ra] and Veech [V1], one analyzes the dynamics of a standard i.e.m.\ $T$ with no connection by considering the first return maps $T^{(n)}$ of $T$ on a decreasing sequence of intervals $I^{(n)}$, with the same left endpoint than $I$. These maps are again standard i.e.m.\ on the same alphabet  $\cA$ but the combinatorial data may be different. The set of all possible combinatorial data accessible from the initial one by this process constitute a  Rauzy class. To each Rauzy class is associated a Rauzy diagram (whose vertices are the elements in the Rauzy class and arrows are the possible transitions). The sequence of combinatorial data for the $T^{(n)}$ is an infinite path in this diagram which can be viewed as a ``rotation number''.

By suspending an i.e.m.\ through Veech zippered rectangle construction [V2], one obtains a linear flow on a translation surface. The genus g of the surface only depends on the Rauzy class.

For a generalized i.e.m.\ $T$ with no connection one can still define the $T^{(n)}$ and obtain an infinite path in a Rauzy diagram. When this path is also associated with a standard i.e.m.\  $T_0$ with  no connection (one then says that $T$ is irrational),  $T$ is semi--conjugate to $T_0$ .

When $d=2$, or more generally $g=1$, such a semi--conjugacy for an affine i.e.m is always a conjugacy as recalled above.

Levitt [L] found an example of an affine irrational i.e.m.\ in higher genus which has a wandering interval. The corresponding standard i.e.m is not unique in his case; this only happens in the non-uniquely ergodic case which has measure zero in parameter space [Ma], [V2].

Later Camelier and Gutierrez [CG] exhibited an example of affine irrational i.e.m.\  with a wandering interval such that the corresponding standard i.e.m.\ is uniquely ergodic. The infinite path in the Rauzy diagram in their case is periodic. The same example was studied more deeply by Cobo [Co]. In particular, he put in evidence on this example the importance of the Oseledets decomposition of the extended Zorich cocycle (see Section 3.1 below). The smoothness of the possible conjugacy between an affine i.e.m.\ and a standard
i.e.m.\ is discussed by Liousse and Marzougui [LM].

Very recently, Bressaud, Hubert and Maass  [BHM] generalized the Camelier-Gutierrez example to a large class of periodic paths in Rauzy diagrams with $g>1$. In the periodic case, the Zorich cocycle is just a matrix in $\hbox{SL}\,(\Z, d)$ with positive coefficients. The vector of the logarithms of the slopes (for the affine i.e.m.\ ) must lie in the Perron-Frobenius hyperplane for this matrix; however, it can have a non-zero component with respect to the next biggest eigenvalue (which is assumed to be real and conjugate to the largest one), and such a choice lead to the required examples.

Our main result is of a similar nature, but instead of starting with periodic paths (a countable set of possibilities) , we consider a set of ``rotation numbers'' of full measure.

Let us fix combinatorial data, such that the associated surface has genus $g>1$. By a deep result of Avila-Viana [AV], the extended Zorich cocycle has $g$ simple positive Lyapunov exponents $\theta_1>\theta_2>\ldots>\theta_g$.  Let $E_0=\R^\cA\supset E_1\supset E_2\supset\ldots\supset E_g$ (with $\dim E_i=d-i$) be the corresponding filtration (defined for almost all parameter values); a necessary and sufficient condition for a vector in $\R^\cA$ to have for coordinates the logarithms of the slopes of an affine i.e.m.\ with this rotation number is that it belongs to the hyperplane $E_1$.

\vskip .3 truecm\noindent {\bf Theorem.}{\it For almost all standard i.e.m. $T_0$ with the given combinatorial data, the following holds: the coordinates of any vector in $E_1\setminus E_2$ can be realized as the logarithms of the slopes of an affine i.e.m.\ semi--conjugate to $T_0$ with a wandering interval. }

\vskip .3 truecm\noindent

We will now summarize the contents of our paper.
In the first section we introduce
interval exchange maps and we develop the continued fraction
algorithms. Accelerating the
Rauzy--Veech map by grouping together arrows with the same type in
the Rauzy diagram leads to the Zorich continued fraction algorithm
(described in 1.2.4) which has the advantage of having a finite
mass a.c.i.m..  The notations and the
presentation of the Rauzy--Veech--Zorich algorithms follow closely
the expository paper [Y1] (see also [Y2]).

Section 2 is devoted to the study of the deformations of affine interval echange maps. First we describe the compact convex set $\hbox{Aff}^{(1)}\, (\underline\gamma ,w)$ of affine i.e.m.\ of the unit interval whose slope vector $w$
and orbit $\underline\gamma$ under the Rauzy--Veech algorithm are prescribed. Following an analogy with the theory of
holomorphic motions in complex dynamics, we them define affine motions. This allows us to characterize the tangent space to $\hbox{Aff}^{(1)}\, (\underline\gamma ,w)$.

In Section 3 deals with the construction of affine interval exchange maps with a wandering interval.

\vskip 1. truecm \noindent {\bf Acknowledgements}  This research
has been supported by the  following institutions: the Coll\`ege de France, the Scuola
Normale Superiore and the Italian MURST (PRIN grant 2007B3RBEY Dynamical Systems and applications) We are also grateful to the two former
institutions and to the Centro di Ricerca Matematica ``Ennio De
Giorgi'' in Pisa for hospitality.

Finally, we thank the referee for his remarks which helped us to improve substantially
the presentation of the paper.

\vfill\eject \noindent {\bf 1. The continued fraction algorithm
for interval exchange maps}

\vskip .5 truecm \noindent {\bf 1.1 Interval exchange maps}

\vskip .3 truecm\noindent  An interval exchange map (i.e.m.\ ) is
determined by combinatorial data on one side, length data on the
other side.

Let $\cal A$ be an alphabet with $d\ge 2$ elements which serve as indices for the intervals.
The combinatorial data is a pair $\pi=(\pi_t,\pi_b)$ of bijections from $\cA$ onto
$\{1, \ldots ,d\}$ which indicates in which order the intervals are met in the domain and in the range of the i.e.m.\ . We always assume that the combinatorial data are {\it irreducible}: for $1\le k<d$, we have
$$
\pi_t^{-1}(\{1,\ldots ,k\})\not=\pi_b^{-1}(\{1,\ldots ,k\})\; .
$$
The length data are the lengths  $(\lambda_\alpha)_{\alpha\in\cA}$ of the subintervals.
Let $T=T_{\pi,\lambda}$ be the i.e.m.\ determined by these data; it is acting on $I=(0,\lambda^*)$, with
$$
\lambda^* = \sum_{\alpha\in\cA}\lambda_\alpha\; .
$$
The subintervals in the domain are
$$
I_\alpha^t = \left( \sum_{\pi_t\beta<\pi_t\alpha}\lambda_\beta , \sum_{\pi_t\beta\le\pi_t\alpha}\lambda_\beta\right)
$$
and those in the range are
$$
I_\alpha^b = \left( \sum_{\pi_b\beta<\pi_b\alpha}\lambda_\beta , \sum_{\pi_b\beta\le\pi_b\alpha}\lambda_\beta\right)
\; .
$$
We also write $I_\alpha$ for $I_\alpha^t$. The translation vector $(\delta_\alpha)_{\alpha\in\cA}$ is given by
$$
\delta_\alpha = \sum_\beta\Omega_{\alpha\beta}\lambda_\beta
$$
where the antisymmetric matrix $\Omega = \Omega (\pi)$ is defined by
$$
\Omega_{\alpha \beta} = \cases{ +1 & if $\pi_t\beta >\pi_t
\alpha \; , \; \pi_b\beta <\pi_b\alpha $, \cr -1 & if
$\pi_t\beta <\pi_t \alpha \; , \; \pi_b\beta >\pi_b \alpha
$, \cr 0 & otherwise. }
$$
We denote the rank of $\Omega$ by $2g$; in fact $g$ is the genus of the translation surfaces obtained from $T$ by suspension. One has thus
$$
\eqalign{
T(x) &= x+\delta_\alpha \; \; \hbox{for}\, x\in I_\alpha^t\, , \cr
T(I_\alpha^t) &= I_\alpha^b\; \; \hbox{for}\, \alpha\in\cA\; . \cr}
$$
We denote by $u_1^t<\ldots <u_{d-1}^t$ the points of $I\setminus\cup_{\alpha\in\cA} I_\alpha^t$, which we call {\it singularities of} $T$. Similarly, the points $u_1^b<\ldots <u_{d-1}^b$ of $I\setminus\cup_{\alpha\in\cA} I_\alpha^b$
are called the {\it singularities of} $T^{-1}$. A {\it connection} is a triple $(u_i^t,u_j^b,m)$, where $m$ is a nonnegative integer, such that
$$
T^m(u_j^b)=u_i^t\; .
$$
Keane has proved [Ke] that an i.e.m.\ with no connection is minimal, and also that an i.e.m.\ has no connection if the length data are independent over $\Q$.

\vskip .5 truecm \noindent {\bf 1.2 The elementary step of the Rauzy--Veech
algorithm}

\vskip .5 truecm \noindent Let $T=T_{\pi, \lambda}$ be an i.e.m.\ . Denote by
$\alpha_t, \alpha_b$ the  elements of
$\cA$ such that
$$
\pi_t(\alpha_t)=\pi_b(\alpha_b)=d\; .
$$
When $u_{d-1}^t\not=u_{d-1}^b$ (which must happen if $T$ has no connection), we consider the first return map $\hat T$
on $\hat I=(0,\hbox{Max}\, (u_{d-1}^t, u_{d-1}^b))$.

When $u_{d-1}^t<u_{d-1}^b$, we have
$$
\hat T(y) = \cases{T^2(y) \; & if $y\in I_{\alpha_b}^t$, \cr
T(y) \; & otherwise.  \cr}
$$
Thus $\hat T$ is an i.e.m.\ with the same alphabet $\cA$, length data $\hat \lambda$, combinatorial data $\hat \pi$ with
$$
\eqalign{
\hat{\lambda}_{\alpha_t} &= \lambda_{\alpha_t}
- \lambda_{\alpha_{b}} \; , \cr
\hat{\lambda}_\alpha & = \lambda_\alpha \;\;,
\alpha\not=\alpha_t\; ,\cr
\hat\pi_t &=\pi_t\; , \cr
\hat\pi_b(\alpha ) &= \cases{
\pi_b(\alpha ) &if $\pi_b(\alpha )\le \pi_b(\alpha_t)$, \cr
\pi_b(\alpha )+1 &if $\pi_b(\alpha_t)<\pi_b(\alpha)<d$, \cr
\pi_b(\alpha_t)+1 &if $\pi_b(\alpha)=d$.\cr}
\cr}
$$

When $u_{d-1}^b<u_{d-1}^t$, we have
$$
\hat T^{-1}(y) = \cases{T^{-2}(y) \; & if $y\in I_{\alpha_t}^b$, \cr
T^{-1}(y) \; & otherwise.  \cr}
$$
In this case, the length and combinatorial data for $\hat T$ are:
$$
\eqalign{
\hat{\lambda}_{\alpha_b} &= \lambda_{\alpha_b}
- \lambda_{\alpha_{t}} \; , \cr
\hat{\lambda}_\alpha & = \lambda_\alpha \;\;,
\alpha\not=\alpha_b\; ,\cr
\hat\pi_b &=\pi_b\; , \cr
\hat\pi_t(\alpha ) &= \cases{
\pi_t(\alpha ) &if $\pi_t(\alpha )\le \pi_t(\alpha_b)$, \cr
\pi_t(\alpha )+1 &if $\pi_t(\alpha_b)<\pi_t(\alpha)<d$, \cr
\pi_t(\alpha_b)+1 &if $\pi_t(\alpha)=d$.\cr}
\cr}
$$
We say that $\hat T$ is deduced from $T$ by an elementary step of the Rauzy--Veech algorithm. We also define
the Rauzy operation $\hat\pi =R_t(\pi)$ (respectively $\hat\pi = R_b(\pi)$) for the change of combinatorial data when $u_{d-1}^t<u_{d-1}^b$ (respectively $u_{d-1}^b<u_{d-1}^t$).

\vskip .5 truecm \noindent {\bf 1.3 Rauzy diagrams }

\vskip .5 truecm \noindent
A {\it Rauzy class} on an alphabet $\cA$ is a nonempty set of irreducible combinatorial data
which is invariant under $R_t,R_b$ and minimal with respect to this property.
A {\it Rauzy diagram} is a graph whose vertices are the elements of a Rauzy class and whose arrows
connect a vertex $\pi$ to its images $R_t(\pi)$ and $R_b(\pi)$. Each vertex is therefore the origin of two arrows. As
$R_t,R_b$ are invertible, each vertex is also the endpoint of two arrows.
It is a fact that the rank of the matrix $\Omega (\pi)$ is the same for all $\pi$ in a given Rauzy class.

An arrow connecting $\pi$ to $R_t(\pi)$ (respectively $R_b(\pi)$) is said to be of {\it top type} (resp.\ {\it bottom type}). The {\it winner} of an arrow of top (resp.\ bottom) type starting at $\pi=(\pi_t,\pi_b)$ with $\pi_t(\alpha_t)=
\pi_b(\alpha_b)=d$ is the letter $\alpha_t$ (resp.\ $\alpha_b$) while the {\it loser} is $\alpha_b$ (resp.\ $\alpha_t$).

To an arrow $\gamma$ of a Rauzy diagram $\cD$ starting at $\pi$ of top (resp.\ bottom) type, is associated
the matrix $B_\gamma\in\hbox{SL}\, (\Z^\cA)$ defined by
$$
B_\gamma = \1 + E_{\alpha_b\alpha_t}
$$
(resp.\ $B_\gamma = \1+E_{\alpha_t\alpha_b}$), where $E_{\alpha\beta}$ is the elementary matrix whose only nonzero coefficient is $1$ in position $\alpha\beta$. For a path $\gamma$ in $\cD$ made of the successive arrows $\gamma_1\ldots \gamma_l$ we associate the product $B_\gamma=B_{\gamma_l}\ldots B_{\gamma_1}$. It belongs to $\hbox{SL}\, (\Z^\cA)$ and has nonnegative coefficients.

A path $\gamma$ in $\cD$ is {\it complete} if each letter in $\cA$ is the winner of at least one arrow in $\gamma$; it is $k$--{\it complete} if $\gamma$ is the concatenation of $k$ complete paths. An infinite path is $\infty$--{\it complete} if it is the concatenation of infinitely many complete paths. By [MMY, Section 1.2.4], if a path $\gamma$ is $(2d-3)$--complete, then all coefficients of $B_\gamma$ are strictly positive.

\vskip .5 truecm\noindent {\bf 1.4 The Rauzy-Veech and Zorich algorithms}

\vskip .5 truecm\noindent Let $T^{(0)}=T_{(\lambda^{(0)}, \pi^{(0)})}$ be an i.e.m.\ with no connection. We
denote by $\cA$ the alphabet for $\pi^{(0)}$ and by $\cD$ the Rauzy diagram on $\cA$ having $\pi^{(0)}$ as a vertex.
The i.e.m.\ $T^{(1)}=T_{(\lambda^{(1)}, \pi^{(1)})}$ deduced from $T^{(0)}$ by the elementary step of the
Rauzy--Veech algorithm has also no connection. It is therefore possible to iterate this elementary step indefinitely
and get a sequence $T^{(n)}=T_{(\lambda^{(n)}, \pi^{(n)})}$ of i.e.m.\ acting on a decreasing sequence $I^{(n)}$ of
intervals and a sequence $\gamma (n,n+1)$ of arrows in $\cD$ from $\pi^{(n)}$ to $\pi^{(n+1)}$. For $m<n$, we also
write $\gamma (m,n)$ for the path from $\pi^{(m)}$ to $\pi^{(n)}$ composed of the $\gamma (l,l+1)$, $m\le l<n$. One has
$$
\eqalign{
\lambda^{(m)} &= ^tB_{\gamma (m,n)}\lambda^{(n)}\, , \cr
\delta^{(n)} &= B_{\gamma (m,n)}\delta^{(m)}\, . \cr}
$$
Conversely, if it is possible to iterate indefinitely the Rauzy--Veech
elementary step starting from $T^{(0)}$, then $T^{(0)}$ has no connection.

Let $\underline\gamma$ be the infinite path starting at $\pi^{(0)}$
obtained by concatenation of the $\gamma (n,n+1)$; then $\underline\gamma$ is
$\infty$--complete. Conversely, if an infinite path $\underline\gamma$ is
$\infty$--complete, it is associated by the Rauzy--Veech algorithm to some
$T=T_{\lambda, \pi}$ with no connection. This $T$ is unique up to rescaling
if and only if it is uniquely ergodic; this last property is true for almost
all $\lambda$ ([Ma], [V2]).

Following Zorich [Z1] it is often convenient to group together in a single
Zorich step successive elementary steps of the Rauzy--Veech algorithm whose
corresponding arrows have the same type (or equivalently the same winner);
we therefore introduce a sequence $0=n_0<n_1<\ldots $ such that for each $k$ all
arrows in $\gamma (n_k,n_{k+1})$ have the same type and this type is alternatively
top and bottom. For $n\ge 0$, the integer $k$ such that $n_k\le n<n_{k+1}$ is called
the {\it Zorich time} and denoted by $Z(n)$.

\vskip .5 truecm\noindent {\bf 1.5 Dynamics of the continued fraction algorithms}

\vskip .5 truecm\noindent
Let $\cR$ be a Rauzy class on an alphabet $\cA$. The elementary step of the Rauzy--Veech algorithm,
$$
(\pi, \lambda) \mapsto (\hat\pi, \hat\lambda )\, ,
$$
considered up to rescaling, defines a map from $\cR\times\P ((\R^+)^\cA)$ to itself, denoted by $Q_{\hbox{\sevenrm RV}}$. There exists a unique absolutely continuous measure invariant under these dynamics ([V2]); it is conservative and ergodic but has infinite total mass, which does not allow all ergodic--theoretic machinery to apply. Replacing a Rauzy--Veech elementary step by a Zorich step gives a new map $Q_{\hbox{\sevenrm Z}}$ on $\cR\times\P ((\R^+)^\cA)$. This map has now a {\it finite} absolutely continuous invariant measure, which is ergodic ([Z1]).

It is also useful to consider the natural extensions of the maps $Q_{\hbox{\sevenrm RV}}$ and $Q_{\hbox{\sevenrm Z}}$, defined through the suspension data which serve to construct translation surfaces from i.e.m.\ . For $\pi\in\cR$, let $\Theta_\pi$ be the convex open cone in $\R^\cA$ defined by the inequalities
$$
\sum_{\pi_t\alpha\le k}\tau_\alpha >0\; , \;\;\; \sum_{\pi_b\alpha\le k}\tau_\alpha <0\; , \;\;\; 1\le k<d\; .
$$
Define also
$$
\eqalign{
\Theta_\pi^t &= \{ \tau\in\Theta_\pi\, , \, \sum_\alpha\tau_\alpha <0\}\; , \cr
\Theta_\pi^b &= \{ \tau\in\Theta_\pi\, , \, \sum_\alpha\tau_\alpha >0\}\; . \cr
}
$$
Let $\gamma\, : \, \pi\rightarrow\hat\pi$ be an arrow in the Rauzy diagram $\cD$ associated to $\cR$. Then $^tB_\gamma^{-1}$ sends $\Theta_\pi$ isomorphically onto $\Theta^t_{\hat\pi}$ (resp.\ $\Theta^b_{\hat\pi}$)
when $\gamma$ is of top type (resp.\ bottom type). The natural extension $\hat Q_{\hbox{\sevenrm RV}}$ is then defined on
$\sqcup_{\pi\in\cR}\{\pi\}\times\P ((\R^+)^\cA)\times\P (\Theta_\pi)$ by
$$
(\pi,\lambda, \tau)\, \mapsto (\hat\pi, ^tB_\gamma^{-1}\lambda , ^tB_\gamma^{-1}\tau )
$$
where $\gamma$ is the arrow starting at $\pi$, associated to the map $Q_{\hbox{\sevenrm RV}}$ at $(\pi,\lambda )$.
The map $\hat Q_{\hbox{\sevenrm RV}}$ has again a unique absolutely continuous invariant measure; it is ergodic, conservative but infinite. One defines similarly a natural extension $\hat Q_{\hbox{\sevenrm Z}}$ for $ Q_{\hbox{\sevenrm Z}}$; it has a unique absolutely continuous invariant measure, which is finite and ergodic.

\vskip .5 truecm\noindent {\bf 1.6 The continued fraction algorithm for generalized and affine i.e.m.\ }

\vskip .5 truecm\noindent
Let $\cA$ be an alphabet and $\pi =(\pi_t,\pi_b)$ be irreducible combinatorial data over $\cA$. Let $I=(0,\lambda^*)$
be an interval and let
$$
\eqalign{
0 &= u_0^t<u_1^t<\ldots <u_d^t=\lambda^*\; , \cr
0 &= u_0^b<u_1^b<\ldots <u_d^b=\lambda^*\; , \cr
}
$$
two sets of points in $\overline I$. Define
$$
\eqalign{
I_\alpha^t &= \left( u^t_{\pi_t(\alpha )-1}, u^t_{\pi_t(\alpha )}\right) \, , \cr
I_\alpha^b &= \left( u^b_{\pi_b(\alpha )-1}, u^b_{\pi_b(\alpha )}\right) \, . \cr
}
$$
A {\it generalized i.e.m.\ } with combinatorial data $\pi$ is a map on $I$ whose restriction to each $I_\alpha^t$ is a
non decreasing homeomorphism onto $I_\alpha^b$ (for some choice of the $u_i^t$, $u_j^b$). When these restrictions are affine, we say that $T$ is an {\it affine} i.e.m.\ .

Connections for generalized i.e.m.\ are again defined by some relation $T^m(u_j^b)$ $=$ $u_i^t$, with $m\ge 0$, $0<i,j<d$. When $T$ has no connection, one has in particular $u^t_{d-1}\not=u^b_{d-1}$. One then defines $\hat I = (0,\hbox{Max}\, (u^t_{d-1},u^b_{d-1}))$ and $\hat T$ as the first return map of $T$ in $\hat I$. Then $\hat T$ is again a generalized i.e.m.\ (affine if $T$ was affine), the combinatorial data being $R_t(\pi)$ if $u^t_{d-1}<u^b_{d-1}$, $R_b(\pi)$ if
$u^b_{d-1}<u^t_{d-1}$. Also, $\hat T$ has no connection, hence we can iterate the processus.

A difference with the case of standard i.e.m. is that the infinite path $\underline\gamma$ in the Rauzy diagram $\cD$ having $\pi$ as a vertex is not always $\infty$--complete.

When this path $\underline\gamma$ is $\infty$--complete, there exists also a standard i.e.m.\ $T_0$ associated to $\underline\gamma$, and any two such $T_0$ are topologically conjugate. Let $I_0$ be the interval on which acts $T_0$. Then there exists a unique semiconjugacy from $T$ to $T_0$, i.e.\ a continuous non--decreasing surjective map $h$ from $I$ onto $I_0$ such that $h\circ T=T_0\circ h$.

\vskip 1. truecm\noindent
{\bf 2. Deformations of affine interval exchange maps}

\vskip .5 truecm\noindent
Let $\cal D$ be a Rauzy diagram on the alphabet $\cA$ and let $\underline\gamma$ be an $\infty$--complete path in $\cal D$ issued from $(\pi_t,\pi_b)$.

An affine i.e.m.\ with combinatorial data $\pi$ is uniquely defined by the lengths $|I^t_\alpha|$ and
$|I^b_\alpha|$ subjected to the only constraint $\sum_\alpha |I^t_\alpha|=\sum_\alpha |I^b_\alpha|$.

Let $w\in\R^\cA$. We will describe the set $\hbox{Aff}\, (\underline\gamma ,w )$ of the affine interval exchange maps whose orbit under the Rauzy--Veech algorithm is given by $\underline\gamma$ and with slope vector $\exp w$:
$$
 |I_\alpha^b|=\exp w_\alpha |I_\alpha^t|\; , \; \forall \alpha\in\cA\; . \leqno{(1)}
$$
We denote by $\hbox{Aff}^{(1)}\, (\underline\gamma ,w)$ the set of affine i.e.m.\ in $\hbox{Aff}\, (\underline\gamma ,w )$ whose domain is $[0,1]$.

When $w=0$ it is known ([Ka], [V1]) that the set of length vectors $\lambda$ corresponding to a fixed Rauzy--Veech expansion $\underline\gamma$ is a simplicial cone of dimension $\le g$ (where $g$ is the genus of the surface associated to the diagram $\cal D$). In the remaining part of Section 2 we assume that $w\not= 0$.

\vskip .5 truecm \noindent
{\bf 2.1 The set $\hbox{Aff}^{(1)}\, (\underline\gamma ,w)$.}

\vskip .3 truecm\noindent
We will first determine a necessary and sufficient condition for $\hbox{Aff}\, (\underline\gamma ,w )\not=\emptyset$.

\vskip .3 truecm \noindent {\bf Lemma 1 } {\it Let $\alpha_t,\alpha_b$ the elements of $\cA$ such that $\pi_t(\alpha_t)=\pi_b(\alpha_b)=d$. There exists an affine interval exchange map of slope $\exp w$ verifying $|I_{\alpha_t}^t|>
|I_{\alpha_b}^b|$ if and only if the intersection}
$$
\{ \sum\lambda_\alpha w_\alpha =0\}\cap\{\lambda_\alpha >0\, , \lambda_{\alpha_t}>\lambda_{\alpha_b} \}
$$
{\it is not empty. }

\vskip .3 truecm\noindent\proof
There exists an affine i.e.m.\ of slope $\exp w$ verifying $|I^t_{\alpha_t}|>|I^b_{\alpha_b}|$ if and only if the hyperplane $\{\sum_\alpha |I^t_\alpha|(\exp w_\alpha -1)=0\}$ intersects the cone $\{|I^t_\alpha|>0\, , \, |I^t_{\alpha_t}|>\exp w_{\alpha_b}|I^t_{\alpha_b}|\}$.

Let $a\not= 0$ in $\R^\cA$. The hyperplane $\{\sum_\alpha a_\alpha x_\alpha =0\}$ does not intersect the positive cone $\{x_\alpha >0\}$ if and only if either all $a_\alpha\ge 0$ or all $a_\alpha\le 0$.

Set first $x_\alpha=|I^t_\alpha|$ for $\alpha\not=\alpha_t$, $x_{\alpha_t}=|I^t_{\alpha_t}|-\exp (w_{\alpha_b})|I^t_{\alpha_b}|$, $a_\alpha = \exp w_\alpha -1$ for $\alpha\not=\alpha_b$, $a_{\alpha_b}=\exp (w_{\alpha_t}+w_{\alpha_b})-1$. We have $\sum_\alpha a_\alpha x_\alpha=\sum_\alpha |I^t_\alpha|(\exp w_\alpha -1)$. Therefore the hyperplane $\{\sum_\alpha |I^t_\alpha|(\exp w_\alpha -1)=0\}$ does not intersect the cone $\{|I^t_\alpha|>0\, , \, |I^t_{\alpha_t}|>\exp w_{\alpha_b}|I^t_{\alpha_b}|\}$ iff
\item{$\bullet$} either $\exp w_\alpha -1\ge 0$ for $\alpha\not=\alpha_b$ and $\exp (w_{\alpha_t}+w_{\alpha_b})-1\ge 0$,
\item{$\bullet$} or $\exp w_\alpha -1\le 0$ for $\alpha\not=\alpha_b$ and $\exp (w_{\alpha_t}+w_{\alpha_b})-1\le 0$.

This is in turn respectively equivalent to
\item{$\bullet$} $w_\alpha\ge 0$ for $\alpha\not=\alpha_b$ and $w_{\alpha_t}+w_{\alpha_b}\ge 0$,
\item{$\bullet$} $w_\alpha\le 0$ for $\alpha\not=\alpha_b$ and $w_{\alpha_t}+w_{\alpha_b}\le 0$.

Take now $x_\alpha=\lambda_\alpha$ for $\alpha\not=\alpha_t$, $x_{\alpha_t}=\lambda_{\alpha_t}-\lambda_{\alpha_b}$; $a_\alpha=w_\alpha$ for $\alpha\not=\alpha_b$,
$a_{\alpha_b}=w_{\alpha_t}+w_{\alpha_b}$. We have $\sum_\alpha a_\alpha x_\alpha=\sum_\alpha\lambda_\alpha w_\alpha$. Therefore the hyperplane $\{\sum_\alpha\lambda_\alpha w_\alpha=0\}$ does not intersect
$\{\lambda_\alpha >0\, , \, \lambda_{\alpha_t}>\lambda_{\alpha_b}\}$ if and only if
\item{$\bullet$} either $w_\alpha \ge 0$ for $\alpha\not=\alpha_b$ and $w_{\alpha_t}+w_{\alpha_b}\ge 0$,
\item{$\bullet$} or $w_\alpha \le 0$ for $\alpha\not=\alpha_b$ and $w_{\alpha_t}+w_{\alpha_b}\le 0$.

We have shown that the negations of both statements considered in the Lemma are equivalent to the same set of inequalities. Hence the proof of the Lemma is complete.
\qed

\vskip .3 truecm\noindent
If an affine interval exchange map verifies (1) and $|I^t_{\alpha_t}|>|I^b_{\alpha_b}|$, one can apply a step of the Rauzy--Veech algorithm. The new affine i.e.m.\ $\hat T$ is the return map of $T$ on $\cup_{\alpha\not=\alpha_b}I^b_\alpha$
and its slope vector $\exp \hat w$ is given by
$$
\eqalign{
\hat w_\alpha &= w_\alpha\; , \; \hbox{if}\, \alpha\not=\alpha_b\, , \cr
\hat w_{\alpha_b} &= w_{\alpha_b}+w_{\alpha_t}\, . \cr}
$$
The corresponding lengths are
$$
\eqalign{
|\hat I^t_\alpha| &= |I^t_\alpha|\; , \; \hbox{if}\, \alpha\not=\alpha_t\, , \cr
|\hat I^t_{\alpha_t}| &= |I^t_{\alpha_t}|-\exp (w_{\alpha_b})|I^t_{\alpha_b}|\, . \cr}
$$
It is easy to check that the maps $\hat T$ obtained in this way (as $T$ varies) are determined by the only constraint
$$
|\hat I^b_\alpha|=\exp \hat w_\alpha |\hat I^t_\alpha|\; . \leqno{(1')}
$$
Moreover, the top Rauzy--Veech operation maps the set
$$
\{ \sum \lambda_\alpha w_\alpha= 0\, , \lambda_\alpha >0\, , \, \lambda_{\alpha_t}>\lambda_{\alpha_b}\}
$$
onto the set
$$
\{ \sum \hat\lambda_\alpha \hat w_\alpha= 0\, , \hat\lambda_\alpha >0\}\; ,
$$
where $\hat\lambda$ is connected to $\lambda$ by the formulas of Section 1.2.

Lemma 1 and the subsequent discussion have a symmetric reformulation for the bottom Rauzy--Veech operation
($|I^t_{\alpha_t}|<|I^b_{\alpha_b}|$, $\lambda_{\alpha_t}<\lambda_{\alpha_b}$).

By applying several times the top or bottom versions of Lemma 1 and the subsequent discussion one  obtains

\vskip .3 truecm \noindent {\bf Lemma 2. } {\it Let $\underline\gamma^*$ be a finite initial segment of $\underline\gamma$.
There exists an affine interval exchange map satisfying (1) whose orbit under the Rauzy--Veech algorithm begins with $\underline\gamma^*$
if and only if the set $\{\sum \lambda_\alpha w_\alpha= 0\, , \lambda_\alpha >0\}$ contains a} standard {\it i.e.m.\ whose expansion under the the Rauzy--Veech algorithm begins with $\underline\gamma^*$.}

\vskip .3 truecm\noindent
We now give a necessary and sufficient condition for $\hbox{Aff}\, (\underline\gamma , w)$ to be non empty.

\vskip .3 truecm \noindent {\bf Proposition  } {\it The set $\hbox{Aff}\, (\underline\gamma ,w)$
is not empty if and only if the hyperplane $\{\sum\lambda_\alpha w_\alpha = 0\}$ contains a standard interval
exchange map
whose Rauzy--Veech expansion is equal to $\underline\gamma$. In this case, the set $\hbox{Aff}^{(1)}\, (\underline\gamma ,w)$, parametrized by the $|I_\alpha^t|$, is convex and compact. }

\vskip .3 truecm\noindent\proof
For $\gamma$ an arrow of $\cD$, we define a matrix $B_\gamma [w]\in\hbox{SL}\,(\R^\cA)$ with nonnegative coefficients in the following way. Let $\pi = (\pi_t,\pi_b)$ be the origin of $\gamma$, $\alpha_t, \alpha_b\in\cA$ such that $\pi_t(\alpha_t)=\pi_b(\alpha_b)=d$. If $\gamma$ is of top type, set
$$
B_\gamma [w] = \1 + \exp w_{\alpha_b} E_{\alpha_b\alpha_t}\; .
$$
If $\gamma$ is of bottom type, set
$$
B_\gamma [w] = \1 + (\exp w_{\alpha_b}-1) E_{\alpha_t\alpha_t}
+ E_{\alpha_t\alpha_b} \; .
$$
Observe that $B_\gamma [0]$ is the matrix $B_\gamma$ introduced in 1.3. The positive coefficients for $B_\gamma [w]$ and $B_\gamma$ appear at the same positions. If $T$ is an affine i.e.m.\ with combinatorial data $\pi$, slope $\exp w$, and $\hat T$ is deduced from $T$ by the Rauzy--Veech operation associated to $\gamma$, the respective lengths $|I^t_\alpha|$, $|\hat I^t_\alpha|$ are related by
$$
|I^t|=^tB_\gamma [w]|\hat I^t|\, ,
$$
in view of the formulas in the discussion following Lemma 1. If $\gamma=\gamma_1\ldots\gamma_l$ is a path in
$\cD$, we define
$$
B_\gamma [w]=B_{\gamma_l} [w_{l-1}]\ldots B_{\gamma_1} [w_{0}]\, ,
$$
with $w_0=w$, $w_j=B_{\gamma_1\ldots\gamma_j}[w]$ for $j>0$. If $\hat T$ is deduced from $T$ by a sequence of Rauzy--Veech operations corresponding to $\gamma$, we still have
$$
|I^t|=^tB_\gamma [w]|\hat I^t|\, .
$$
Observe also that the positive coefficients of $B_\gamma [w]$ and $B_\gamma [0]=B_\gamma$ appear again at the same positions. Let now be $\underline\gamma$ be an $\infty$--complete path in $\cD$. Let $\hbox{Aff}\, (\underline\gamma (0,n),w)$ be the set of lengths $(I^t_\alpha)$ for affine i.e.m.\ $T$ whose Rauzy--Veech expansion starts with the initial segment $\underline\gamma (0,n)$ of $\underline\gamma$. We have
$$
\eqalign{
\hbox{Aff}\, (\underline\gamma (0,n),w) &= ^tB_{\underline\gamma (0,n)}[w]((\R^+)^\cA)\, , \cr
\hbox{Aff}\, (\underline\gamma (0,n),0) &= ^tB_{\underline\gamma (0,n)}((\R^+)^\cA)\, , \cr
\hbox{Aff}\, (\underline\gamma ,w) &= \cap_{n\ge 0} \hbox{Aff}\, (\underline\gamma (0,n),w) \, , \cr
\hbox{Aff}\, (\underline\gamma ,0) &= \cap_{n\ge 0} \hbox{Aff}\, (\underline\gamma (0,n),0) \, . \cr
}
$$
Let $n>m$ be such that $\gamma (m,n)$ is $(2d-3)$--complete. Then, as recalled in Section 1.3, all coefficients of $B_{\gamma (m,n)}$ are positive. Therefore, the same is true for $B_{\gamma (m,n)}[
B_{\gamma (0,m)}w]$. We therefore have
$$
\eqalign{
\overline{\hbox{Aff}\, (\underline\gamma (0,n),0)}&=^tB_{\gamma (0,n)}(\overline{(\R^+)^\cA}) \cr
&=
^tB_{\gamma (0,m)}(\overline{^tB_{\gamma (m,n)}((\R^+)^\cA)}\subset \{0\}\cup \hbox{Aff}\, (\gamma (0,m),0)\, , \cr
}
$$
and similarly
$$
\overline{\hbox{Aff}\, (\underline\gamma (0,n),w)}\subset \{0\}\cup \hbox{Aff}\, (\gamma (0,m),w)\; .
$$
It follows that
$$
\eqalignno{
\{0\}\cup \hbox{Aff}\, (\underline\gamma ,0) &= \cap_{n\ge 0} \overline{\hbox{Aff}\, (\underline\gamma (0,n),0)}\, , &(2)\cr
\{0\}\cup \hbox{Aff}\, (\underline\gamma ,w) &= \cap_{n\ge 0} \overline{\hbox{Aff}\, (\underline\gamma (0,n),w)}\, . &(3)\cr}
$$
We conclude that $\hbox{Aff}\, (\underline\gamma ,w)$ is nonempty if and only if $\hbox{Aff}\, (\underline\gamma (0,n),w)$ is nonempty for all $n\ge 0$; by Lemma 2 this happens if and only if $\hbox{Aff}\, (\underline\gamma (0,n),0)$ intersects the hyperplane $\{\sum_\alpha\lambda_\alpha w_\alpha =0\} $ for
all $n\ge 0$; in view of the formula (2) above, this last condition holds if and only if the hyperplane
$\{\sum_\alpha\lambda_\alpha w_\alpha =0\} $ meets $\hbox{Aff}\, (\underline\gamma ,0)$. This proves the first statement in the proposition.

The second statement follows from formula (3) and the fact that
$$
\overline{\hbox{Aff}\, (\underline\gamma (0,n),w)}=^tB_{\gamma (0,n)}[w](\overline{(\R^+)^\cA})
$$
is a closed convex cone for $n\ge 0$. \qed

\vskip .3 truecm\noindent
When there exists a unique (up to rescaling) standard i.e.m.\ whose expansion under the Rauzy--Veech algorithm is $\underline\gamma$ the condition stated in the Proposition above means that the vector $w$ belongs to the hyperplane
$$
\{\sum\lambda_\alpha w_\alpha = 0\}\; .
$$
In general, as already mentioned, $\hbox{Aff}\, (\underline\gamma ,0)$ is a simplicial cone of dimension
$r\le g$. Let us denote by $\lambda^{(1)},\ldots ,\lambda^{(r)}$  the normalized extremal vectors of this simplicial cone. The necessary and sufficient condition which guarantees that $\hbox{Aff}\, (\underline\gamma , w )$ is not empty is that the numbers
$$
\sum_{\alpha\in\cA}\lambda_\alpha^{(j)} w_\alpha\; , \;\; j=1,\ldots , r
$$
are neither all strictly positive, nor all strictly negative.

\vskip .3 truecm\noindent {\bf Remark.} For fixed combinatorial data, normalized affine i.e.m.\ form a manifold of dimension $(2d-2)$, and the standard i.e.m.\ have dimension $(d-1)$. As almost all i.e.m.\
are uniquely ergodic, one can think that $(d-1)$ is also the ``dimension'' of the set of paths $\underline\gamma$. When $\underline\gamma$ corresponds to a uniquely ergodic standard i.e.m.\ , the constraint $\sum_\alpha\lambda_\alpha w_\alpha=0$ defines a $(d-1)$ dimensional space. Therefore one can expect that for most $(\underline\gamma , w)$ the set $\hbox{Aff}^{(1)}\, (\underline\gamma ,w)$ is of dimension $(2d-2)-(d-1)-(d-1)=0$. As $\hbox{Aff}^{(1)}\, (\underline\gamma ,w)$ is convex and compact this would mean that $\hbox{Aff}^{(1)}\, (\underline\gamma ,w)$ is reduced to a point. The problem with this heuristic argument is that the map which associates to an affine i.e.m.\ $T$ its ``rotation number''
$\underline\gamma$ is not smooth.

\vskip .5 truecm \noindent
{\bf 2.2 Affine motions.}

\vskip .3 truecm\noindent
Let $w\not= 0$ and $T^*\in\hbox{Aff}^{(1)}\, (\underline\gamma , 0)$ such that
$$
\sum_\alpha\lambda^*_\alpha w_\alpha =0\; .
$$
We choose an affine i.e.m\ $T_0$ in the intrinsic interior of the nonempty compact
convex set $\hbox{Aff}^{(1)}\, (\underline\gamma , w)$.

We denote by $u_i^b,u_i^t$ ($1\le i\le d-1$) the singularities of $T_0^{-1}$ and $T_0$ respectively.
Let
$$
Z=\{u^t_{i,n},u^b_{j,m}\, ,\, n\le 0\, ,\, m\ge 0\, , \, 1\le i,j\le d-1\}\cup\{0,1\}\,
$$
where we have set $u^t_{i,n}=T_0^n(u_i^t)$, $u^b_{j,m}=T_0^m(u^b_j)$.
As $T^*$ is minimal, the complement in $[0,1]$ of the closure $\overline{Z}$ is the union
of the interiors of the wandering intervals for $T_0$.

Let ${\cal T}(T_0)$ be the space of functions in $L^\infty ([0,1])$ which are $T_0$--invariant, and constant on each wandering interval of $T_0$. We will show below that ${\cal T}(T_0)$ is finite--dimensional and that
there exists a canonical affine isomorphism between $\hbox{Aff}^{(1)}\, (\underline\gamma , w )$ and
a simplex in the hyperplane ${\cal T}_0(T_0)$ of functions in ${\cal T}(T_0)$ with mean value $0$.

Let $T_1\in\hbox{Aff}^{(1)}\, (\underline\gamma , w )$, and let $(T_s)_{s\in [0,1]}$ be the affine
segment (with affine parametrization) joining $T_0$ and $T_1$ in $\hbox{Aff}^{(1)}\, (\underline\gamma , w )$.

Let $u_i^t(s)$ and $u_i^b(s)$ denote the singularities of $T_s$ and $T_s^{-1}$ respectively.
For $n\le 0$, $m\ge 0$, set
$$
\eqalign{
u_{i,n}^t(s) &= T_s^n (u_i^t(s)) \, , \cr
u_{j,m}^b(s) &= T_s^m (u_j^b(s)) \, . \cr
}
$$
As the parametrization of the segment is affine, there exists a function $\nu$ on $Z$ with $\nu (0)=\nu (1)= 0$ such that
$$
\eqalign{
u_{i,n}^t(s) &= u_{i,n}^t+s\nu (u_{i,n}^t)\, , \; \cr
u_{j,m}^b(s) &= u_{j,m}^b+s\nu (u_{j,m}^b)\, , \; \cr
}
$$
for $1\le i,j\le d-1$, $0\le s\le 1$, $n\le 0$, $m\ge 0$.

\vskip .5 truecm\noindent {\bf Lemma.}{\it The function  $\nu$ is Lipschitz on $Z$
and satisfies $x_1+\nu (x_1)>x_0+\nu (x_0)$ for $x_1>x_0$ in $Z$. }

\vskip .3 truecm\noindent\proof
First observe that, as $T_0$ belongs to the interior of $\hbox{Aff}^{(1)}\, (\underline\gamma , w )$,
we can extend the segment $T_s$ to $s\in [a,1]$ for some $a<0$. As all $T_s$ are semiconjugate to
$T^*$, for any $x_0<x_1$ in $Z$, we must have $x_0+s\nu (x_0)<x_1+s\nu (x_1)$ for all $s\in [a,1]$, and
the assertions of the Lemma follow. \qed

\vskip .3 truecm\noindent
We first extend $\nu$ on $\overline{Z}$ by continuity and them (if $\overline{Z}\not=[0,1]$) to the
full interval $[0,1]$ by forcing $\nu$ to be affine on each wandering interval. Let $\mu$ be the
distributional derivative of $\nu$.

\vskip .5 truecm\noindent {\bf Proposition 1.}{\it The map $T_1\mapsto \mu$ is an affine isomorphism
of $\hbox{Aff}^{(1)}\, (\underline\gamma , w )$ onto the subset $\{\mu\in{\cal T}_0(T_0),\,
\mu\ge -1\}$.}

\vskip .3 truecm\noindent\proof
In view of the Lemma, the function $\nu$ constructed above is Lipschitz on $[0,1]$ and satisfies
$x_1+\nu (x_1)\ge x_0+\nu (x_0)$ for $0\le x_0\le x_1\le 1$. Therefore $\mu\in L^\infty([0,1])$
and $\mu\ge -1$. It is also clear that $\mu$ is constant on each wandering interval and has mean
value $0$ (because $\nu (0)=\nu (1)=0$). For $0\le s<1$, define
$$
h_s(x)=x+s\nu (x)\, .
$$
Then $h_s$ is an homeomorphism of $[0,1]$, and from the construction of $\nu$ we have
$$
h_s\circ T_0=T_s\circ h_s\, .
$$
Write $\chi_\alpha$ for the constant value of ${\partial\over\partial s}T_s\vert_{s=0}$
on $I_\alpha^t(T_0)$. Taking the derivative w.r.t.\ $s$ of the last equality, we obtain
$$
\nu(T_0(x))=\nu (x)\exp w_\alpha+\chi_\alpha\, , x\in I_\alpha^t(T_0)\, .
$$
Taking now the derivative w.r.t.\ $x$, we see that $\mu$ is $T_0$--invariant. We have shown so
far that the map $T_1\mapsto\mu$ takes indeed its values in $\{\mu\in {\cal T}_0(T_0), \,
\mu\ge -1\}$.

Conversely, let $\mu\in {\cal T}_0(T_0)$ with $\mu\ge -1$.
Let $\nu$ be the primitive of $\mu$ which vanishes at $0$ and $1$. For $0\le s\le 1$,
define
$$
h_s(x)=x+s\nu (x)\, .
$$
For $0\le s<1$, $h_s$ is an homeomorphism of $[0,1]$. The map $h_1$ is non--decreasing, Lipschitz and
satisfies $h_1(0)=0$, $h_1(1)=1$. For $0\le s<1$, set
$$
T_s=h_s\circ T_0\circ h_s^{-1}\, .
$$
This defines a generalized i.e.m.\ . The $T_0$--invariance of $Dh_s=1+s\mu$ implies that $T_s$ is in fact
affine and lies in $\hbox{Aff}^{(1)}\, (\underline\gamma , w )$. Moreover, the map $s\mapsto T_s$ is affine
on $[0,1)$ and therefore can be extended at $s=1$ to define $T_1\in \hbox{Aff}^{(1)}\, (\underline\gamma , w )$. One has $h_1\circ T_0=T_1\circ h_1$.

Finally, it is clear that the map $\mu\mapsto T_1$ is inverse to the one considered above, and that these
maps are affine. \qed

\vskip .3 truecm \noindent
In the end of this section, we give a description of the space ${\cal T}(T_0)$ and of the subset
$\{\mu\in{\cal T}_0(T_0),\,
\mu\ge -1\}$ appearing in Proposition 1.

It follows from Proposition 1 that ${\cal T}(T_0)$ is finite--dimensional; moreover the dimension
of ${\cal T}_0(T_0)$ is the same than the dimension $r(\underline\gamma , w )-1$ of the affine
subspace supporting $\hbox{Aff}^{(1)}\, (\underline\gamma , w )$. We therefore have
$$
\hbox{dim}\, {\cal T}(T_0)=r(\underline\gamma , w )\le d\, .
$$
Associated to the partition of $[0,1]\, (\hbox{mod}\, 0)$ between $\overline Z$ and $[0,1]\setminus
\overline Z$, we have a direct sum decomposition
$$
{\cal T}(T_0)={\cal T}_c(T_0)\oplus{\cal T}_d(T_0)\,
$$
where:
\item{$\bullet$} ${\cal T}_c(T_0)$ is the subspace of ${\cal T}(T_0)$
formed by functions supported on $\overline Z$; it is non zero if and only if
$\hbox{Leb}\, (\overline Z)>0$;
\item{$\bullet$} ${\cal T}_d(T_0)$ is the subspace of ${\cal T}(T_0)$
formed by functions vanishing on $\overline Z$; it is non zero if and only if
$\overline{ Z}\not=[0,1]$, i.e.\ $T_0$ has wandering intervals.

Denote by $r_d$ (resp.\ $r_c$) the dimension of ${\cal T}_d(T_0)$ (resp.\
${\cal T}_c(T_0)$). We have $r_d+r_c=r(\underline \gamma ,w)\le d$. Obviously,
as a function in ${\cal T}_d(T_0)$ is constant on each wandering interval
and $T_0$--invariant, the dimension $r_d$ is exactly the number of orbits of maximal
wandering intervals. We conclude in particular that this number is finite and $\le d$.

As for the dimension $r_c$ of ${\cal T}_c(T_0)$, observe that a function on $\overline Z$ is $T_0$--invariant if and only if it is measurable w.r.t.\ the $\sigma$--algebra of $T_0$--invariant subsets of
$\overline Z$. (Functions and subsets are considered modulo subsets of Lebesgue measure $0$; this makes
sense as this class of subsets is preserved by $T_0$). Therefore this $\sigma$--algebra is generated
by atoms $Z_1$, $Z_2$, $Z_{r_c}$ which form a partition $\hbox{mod}\, 0$ of $\overline Z$.

It is now easy to describe the subset $\{\mu\in {\cal T}_0(T_0)\, \, \mu\ge -1\}$ of Proposition 1.
It is the simplex whose vertices $\mu_1,\ldots ,\mu_r$ ($r=r(\underline \gamma ,w)$) are defined by
$$
\mu_i(x) = \cases{
-1 &for $x\notin Z_i$, \cr
|Z_i|^{-1}\sum_{j\not=i}|Z_j| &for $x\in Z_i$, \cr}
$$
where the $Z_i$ for $r_c<i\le r$ denote the orbits of the maximal wandering intervals. We have thus proved

\vskip .5 truecm\noindent {\bf Proposition 2.}{\it The set $\hbox{Aff}^{(1)}\, (\underline\gamma , w )$
is a simplex whose vertices are canonically associated to either the orbits of the maximal wandering
intervals of $T_0$ or the atoms amongst the $\hbox{mod}\, 0\; T_0$--invariant subsets of $\overline Z$.}

\vskip .3 truecm\noindent
The correspondence is explained in Proposition 1: if $T$ is a vertex of $\hbox{Aff}^{(1)}\, (\underline\gamma , w )$, there is a semiconjugacy $h_T$ from $T_0$ to $T$ such that $h_T(Z_i)$ has measure $0$ except for
one value of $i$ where it has full measure.

More generally, Proposition 1 defines for each $T\in \hbox{Aff}^{(1)}\, (\underline\gamma , w )$ a
semiconjugacy $h_T$ from $T_0$ to $T$ (which is a conjugacy when $T$ belongs to the interior of
$\hbox{Aff}^{(1)}\, (\underline\gamma , w )$).

On the other hand, for any fixed $x\in [0,1]$, the map $T\mapsto h_T(x)$ is affine. In analogy with the
notion of holomorphic motion, one can speak of an {\it affine motion} of $[0,1]$ parametrized by
$\hbox{Aff}^{(1)}\, (\underline\gamma , w )$. The functions $\mu$ then play the role of
``Beltrami forms''.

\vskip .3 truecm\noindent {\bf Remark.}
In the more familiar case where $w=0$, if one takes $T_0$ in the intrinsic interior of
$\hbox{Aff}^{(1)}\, (\underline\gamma , 0 )$, we have obviously ${\cal T}_d(T_0)=0$ as $T_0$ is minimal.
Lebesgue measure $m$ can be written as a barycentric combination $m=\sum_{i=1}^rt_i\mu_i$, where
$\mu_1,\ldots ,\mu_r$ are the ergodic probability measures invariant under $T_0$. The coefficients $t_i$ are
positive because $T_0$ belongs to the interior of $\hbox{Aff}^{(1)}\, (\underline\gamma , 0 )$.
Each $\mu_i$ is absolutely continuous w.r.t.\ $m$, and they are mutually singular; therefore they are
supported by subsets $Z_1,\ldots ,Z_r$ which form a partition of $[0,1] \, (\hbox{mod}\, 0)$.

In the general case, Lebesgue measure is no more $T_0$--invariant but it is still quasi--invariant.

\vskip 1. truecm\noindent
\noindent {\bf 3. Wandering intervals for affine interval exchange maps}

 \vskip .5 truecm\noindent
 {\bf 3.1 The Zorich cocycle}

\vskip .5 truecm\noindent
Let $\cR$ be a Rauzy class on an alphabet $\cA$, $\cD$ the associated Rauzy diagram. For $T=T_{\pi,\lambda }$
a standard i.e.m.\ acting on some interval $I$ with combinatorial data $\pi\in\cR$, define $E_T$ to be the vector space of functions on $I$ which are constant on each subinterval $I_\alpha^t$. This vector space is canonically isomorphic to
$\R^\cA$. Let $\hat T=T_{\hat \pi ,\hat \lambda }$ be the i.e.m.\ deduced from $T$ by one step of the Rauzy--Veech algorithm, let $\gamma$ be the corresponding arrow from $\pi$ to $\hat\pi$ in $\cD$, let $\hat I$ be the interval on which $\hat T$ acts  and $\hat I_\alpha^t$ the associated subintervals. For $\varphi\in E_T$, one defines a function $\hat\varphi\in E_{\hat T}$ by
$$
\hat \varphi (x) = \sum_{i=0}^{q(x)-1}\varphi (T^ix)\; ,
$$
where $q(x)$ is the return time of $x$ in $\hat I$ (equal to $1$ or $2$). The matrix of the linear map $\varphi\mapsto\hat\varphi$ from $E_T$ to $E_{\hat T}$ in the canonical bases of these spaces is $B_\gamma$.

At the projective level, the fibered map
$$
\eqalign{
(\pi ,\lambda ,\varphi )\, & \mapsto (Q_{\hbox{\sevenrm RV}} (\pi ,\lambda),B_\gamma\varphi )\, , \cr
\cR\times \P ((\R^+)^\cA) \times\R^\cA & \rightarrow \cR\times \P ((\R^+)^\cA) \times\R^\cA\, , \cr}
$$
is called the {\it extended Zorich cocycle} over the Rauzy--Veech dynamics $Q_{\hbox{\sevenrm RV}}$.

There is an invariant subbundle under this cocycle whose fiber over $(\pi ,\lambda)$ is
$$
H(\pi ) = \hbox{Im}\, \Omega (\pi)\; .
$$
Indeed, we have
$$
B_\gamma\Omega (\pi)=\Omega (\hat\pi)\,  ^tB_\gamma^{-1}\, .
$$
It also follows that the restriction of the cocycle to this subbundle, called the {\it Zorich cocycle}, is
symplectic (for the symplectic form defined by the $\Omega (\pi )$). To analyze the extended Zorich cocycle,
one goes to the accelerated dynamics $Q_{\hbox{\sevenrm Z}}$, i.e.\ one reparametrizes the time in the algorithm in order to apply the Oseledets multiplicative ergodic theorem. Then, the Lyapunov exponents on the quotient $\R^\cA/H(\pi )$ are all equal to zero. Avila--Viana ([AV], see also [Fo]) have proved that the Lyapunov exponents on $H(\pi )$ are all simple, hence by symplecticity they can be written as
$$
\theta_1>\theta_2>\ldots >\theta_g>-\theta_g>\ldots >-\theta_1\; .
$$
Here $g={1\over 2}\hbox{dim}\, H(\pi )$ is the genus of the surface obtained by suspension. Associated to these exponents, we have for almost all $T$ a filtration
$$
E_T=\R^\cA=E_0\supset E_1\supset\ldots\supset E_g\, ,
$$
with $\hbox{dim}\, E_i=d-i$. Here, we have
$$
E_1 = \{\varphi\in E_T\, , \,  \int_I\varphi (x)dx=0\}\, .
$$

\vskip .5 truecm \noindent {\bf 3.2 Statement of the main result}

\vskip .3 truecm\noindent
We assume $g\ge 2$. We recall the statement of the Theorem in the introduction.

\vskip .3 truecm\noindent
{\bf Theorem.}{\it For all vertices $\pi$ of $\cal D$, for almost all $\lambda\in (\R^+)^\cA$, for any $w\in E_1(\pi,\lambda )\setminus E_2(\pi,\lambda )$, there exists an affine i.e.m.\ $T^* =T^*_{\pi,\lambda,w}$ with the following properties:
\item{(i)} $T^*\in \hbox{Aff}\, (\underline\gamma , w)$;
\item{(ii)} $T^*$ has a wandering interval.}

\vskip .5 truecm\noindent
{\bf Remarks.} \item{$1.$} For almost all $(\pi,\lambda )$, $T_{\pi,\lambda}$ is uniquely ergodic; then $w\in E_1(\pi,\lambda )$ is a necessary condition
for an affine i.e.m.\ to satisfy (i).
\item{$2.$} Actually the proof of the theorem shows that any affine i.e.m.\ in $\hbox{Aff}\, (\underline\gamma , w)$ has a wandering interval: see the remark at the end of Section 3.7. Moreover, in view of this
remark and of the remark at the end of Section 2.1, it appears very probable that there is up to scaling only one affine i.e.m.\ in $\hbox{Aff}\, (\underline\gamma , w)$

\vskip .5 truecm \noindent {\bf 3.3 Reduction to a statement on Birkhoff sums}

\vskip .5 truecm \noindent {\bf 3.3.1 } The main step in the proof of the theorem is the following result

\vskip .3 truecm\noindent {\bf Proposition.}{\it For all vertices $\pi$ of $\cal D$, for almost all $\lambda\in (\R^+)^\cA$, for all $w\in E_1(\pi,\lambda )\setminus E_2(\pi,\lambda )$, there exists $x^*$, not in the orbits of the singularities of $T_{\pi,\lambda}^{\pm 1}$, such that the Birkhoff sums of $w$ at
$x^*$ satisfy, for all $\varepsilon >0$ and a constant $C(\varepsilon )>0$ independent of $n\in\Z$,}
$$
S_n w(x^*)\le C(\varepsilon )-|n|^{\theta_2/\theta_1 -\varepsilon}\; . $$

\vskip .5 truecm \noindent
The Birkhoff sums are here defined as usual as
$$ S_nw (x^*) = \cases{ \sum_{i=0}^{n-1} w_{\beta_i} & for $n\ge 0$, \cr
                        - \sum_{i=n}^{-1} w_{\beta_i} & for $n< 0$, \cr} $$
with $T_{\pi,\lambda}^i(x^*)\in I^{t}_{\beta_i}$.

\vskip .5 truecm \noindent {\bf 3.3.2 } The theorem follows from the proposition by the usual Denjoy construction. Let
$\pi,\lambda,w,x^*$ be as in the Proposition and $I^{(0)}$ be the interval of definition of $T_{\pi ,\lambda}$. Define, for $n\in\Z$
$$
l_n = \exp\{ S_nw(x^*)\} \, .
$$
From the Proposition it follows that
$$
L=\sum_{n\in\Z} l_n<+\infty\, .
$$
For $x\in I^{(0)}$ set
$$\eqalign{
l^-(x) &= \sum_{T_{\pi,\lambda}^n(x^*)<x} l_n\, , \cr
l^+(x) &= \sum_{T_{\pi,\lambda}^n(x^*)\le x} l_n\, , \cr} $$
and let $h\, : \, [0,L]\rightarrow I^{(0)}$ be the continuous non decreasing map such that
$$ h^{-1}(x) =[l^-(x),l^+(x)]\, . $$
One then defines the affine i.e.m.\ $T^*$ on $[0,L]$ by
\item{$\bullet$} $T^*(l^\pm (x))=l^\pm (T_{\pi,\lambda}(x))$,
\item{$\bullet$} when $l^-(x)<l^+(x)$, $T^*$ is affine from the interval $[l^-(x),l^+(x)]$ onto the interval
$[l^-(T_{\pi,\lambda}(x))$,$l^+(T_{\pi,\lambda}(x))]$.

\noindent Then, the fact that $T^*$ is an affine i.e.m.\ with the required slopes follow from the definition of the $l_i$. The semi--conjugacy to
$T_{\pi,\lambda}$ is built in the construction (using also that $T_{\pi,\lambda}$ is minimal). Finally, the interval $h^{-1}(x^*)$ is wandering.

\vskip .5 truecm \noindent {\bf 3.4 Limit shapes for Birkhoff sums}

\vskip .5 truecm \noindent {\bf 3.4.1 } In order to prove the Proposition in 3.3.1, we construct some functions closely related to the Zorich cocycle. Such functions have also been considered in a different setting in [BHM]. Instead of acting on $(\pi, \lambda)$ we consider the natural extension of the Rauzy--Veech dynamics (and the Zorich acceleration) acting on $(\pi,\lambda,\tau)$, where $\tau\in\R^\cA$ is a suspension datum satisfying the usual conditions (for $1\le k\le d$)
$$\sum_{\pi_t\alpha <k}\tau_\alpha  >0\; , \;\;\;\sum_{\pi_b\alpha <k}\tau_\alpha  <0\; . $$
Instead of a filtration
$$
E_0=\R^\cA\supset E_1(\pi,\lambda)\supset E_2(\pi,\lambda)\supset\ldots
$$
as above, we get from Oseledets theorem $1$--dimensional subspaces $F_i(\pi,\lambda,\tau)$ associated to the Lyapunov exponent $\theta_i$, generated by a vector in $E_{i-1}(\pi,\lambda)\setminus E_i(\pi,\lambda)$. Moreover the sums
$\bigoplus_{j=1}^i F_j(\pi,\lambda,\tau)$ {\it depend only on} $(\pi,\tau)$. (This is the subspace of vectors decreasing in the past  under the Zorich cocycle at a rate at least $-\theta_i$).

In particular $F_1$ depends only on $(\pi,\tau)$, not on $\lambda$; because the matrices $B$ of the Zorich cocycle only have non negative entries (and positive entries after appropriate iteration), the subspace $F_1(\pi,\tau)$ is contained in the positive cone $(\R^+)^\cA$; we write $q(\pi,\lambda)$ for a {\it positive} vector generating $F_1(\pi,\tau)$, normalized by
$$
\sum_\alpha q_\alpha^2(\pi,\tau) =1\, .
$$
Next, we consider the $2$--dimensional subspace $F_1\oplus F_2$, depending only on $(\pi,\tau)$: we choose a vector $v(\pi,\tau)$ satisfying
$$\eqalign{
\sum_\alpha v_\alpha^2(\pi,\tau ) &=1\, , \cr
\sum_\alpha v_\alpha (\pi,\tau)q_\alpha (\pi,\tau) &=0\, . \cr}
$$
There are two choices for $v$, differing by a sign, both of them being relevant in the following; we fix such a choice.

From $q$ and $v$, it is easy to find a generator $w$ for $F_2(\pi,\lambda,\tau)$. Indeed we have
$$
F_2(\pi,\lambda,\tau)\subset E_1 (\pi,\lambda)\, , $$
with
$$
E_1 (\pi,\lambda) = \{ w\, , \, \sum_\alpha\lambda_\alpha w_\alpha =0\}\, . $$
Therefore, we will take
$$
w(\pi,\lambda,\tau)=v(\pi,\tau)-t(\pi,\lambda,\tau)q(\pi,\tau)$$
with
$$
t(\pi,\lambda,\tau)={\langle\lambda,v\rangle\over\langle\lambda, q\rangle}\, .
$$

\vskip .3 truecm\noindent {\bf Proposition.}{\it $\,$ For almost all $(\pi,\lambda,\tau)$ and all $(n_\alpha)\in\N^\cA$, not all equal to $0$, we have}
$$
\sum_\alpha n_\alpha w_\alpha(\pi,\lambda,\tau)\not= 0\, .
$$

\vskip .5 truecm \noindent\proof Indeed, fixing $(n_\alpha)$, we have
$$
\sum_\alpha n_\alpha w_\alpha =0  \Leftrightarrow t = {\langle n,v\rangle\over\langle n, q\rangle}\,
$$
where $\langle n, q\rangle >0$ as $n_\alpha\ge 0$, $q_\alpha >0$. In view of the formula for $t$, for
fixed $(\pi ,\lambda)$ this happens with measure $0$ w.r.t.\ $\lambda$. The conclusion follows by Fubini's theorem. \qed

\vskip .5 truecm \noindent {\bf 3.4.2 The functions $V_\alpha (\pi ,\lambda)$.} Let $(\pi,\tau)$ be a typical point (for backward time Rauzy--Veech--Zorich dynamics). Let $(\pi^{(-n)},\tau^{(-n)})$ be its backwards orbit for the Rauzy--Veech dynamics. Let $q^{(-n)}(\pi,\tau )$, $v^{(-n)}(\pi,\tau )$ be the images of $q(\pi,\tau )$, $v(\pi,\tau )$
under the Zorich cocycle. From the invariance of $F_1$ and $F_1\oplus F_2$ w.r.t.\ the Zorich cocycle we can write
$$\eqalign{
q^{(-n)}(\pi,\tau ) &= \Theta_1^{(-n)} q(\pi^{(-n)},\tau^{(-n)} )\, , \cr
v^{(-n)}(\pi,\tau ) &= \Theta_2^{(-n)} v(\pi^{(-n)},\tau^{(-n)} ) + \Theta^{(-n)} q(\pi^{(-n)},\tau^{(-n)} ) \, , \cr}
$$
where $\Theta_1^{(-n)}$, $\Theta_2^{(-n)}$ and $\Theta^{(-n)}$ are real numbers depending on $\pi ,\tau ,n$,
$\Theta_1^{(-n)}>0$. We will always make a coherent choice for the vectors $v(\pi^{(-n)},\tau^{(-n)} )$ along an orbit in order to have $\Theta_2^{(-n)}>0$.
The coefficient
$\Theta_1^{(-n)}$ is exponentially small (in Zorich reparametrized time) at rate $\theta_1$,
$\Theta_2^{(-n)}$ is exponentially small at rate $\theta_2$, and $|\Theta^{(-n)}|$ is at most exponentially small at rate $\theta_2$.

Let $u^{(-n)}(\pi,\tau )=(q^{(-n)}(\pi,\tau ),v^{(-n)}(\pi,\tau ))$. According to the definition of the Zorich cocycle, we have
$$
u^{(-n)}_\beta = u^{(-n-1)}_\beta \, ,
$$
if $\beta$ is {\it not} the loser of the arrow from $\pi^{(-n-1)}$ to $\pi^{(-n)}$ and
$$
u^{(-n)}_{\beta_l} = u^{(-n-1)}_{\beta_l}+u^{(-n-1)}_{\beta_w}\, ,
$$
if $\beta_l$ (resp.\ $\beta_w$) is the loser (resp.\ the winner) of this arrow.

For $\alpha\in\cA$, let $\Gamma_\alpha^{(-n)}$ be the broken line in $\R^2$ starting at the origin and obtained by adding successively the vectors $ u^{(-n)}_{\beta_i}$, where $\beta_0, \beta_1, \ldots$ are defined as follows: if $T^{(0)}$ is any i.e.m.\ with combinatorial data $\pi^{(0)}$, and $T^{(-n)}$ is the i.e.m.\ whose $n$--times Rauzy--Veech induction is $T^{(0)}$, we have
$$
[T^{(-n)}]^i(I^{(0)}_\alpha )\subseteq I^{(-n)}_{\beta_i}\; . $$
Here, $i$ runs from $0$ to the return time of $I_\alpha^{(0)}$ in $I^{(0)}$.

In other terms, $\beta_0, \beta_1, \ldots$ is the itinerary of $I_\alpha^{(0)}$ with respect to the partition of $I^{(-n)}$ by the $I_\beta^{(-n)}$. When we go one step further to $T^{(-n-1)}$ on $I^{(-n-1)}$, the new itinerary is obtained by replacing $\beta_l$ by $\beta_l\beta_w$ or $\beta_w\beta_l$ (depending whether the arrow from $\pi^{(-n-1)}$ to $\pi^{(-n)}$ has top or bottom type).

Consequently, the vertices of $\Gamma_\alpha^{(-n)}$ are also vertices of $\Gamma_\alpha^{(-n-1)}$.
The following properties are now clear:
\item{$1$.} $\Gamma_\alpha^{(-n)}$ is the graph of a piecewise affine continuous map $V_\alpha^{(-n)}(\pi,\tau)$ on $[0,q_\alpha(\pi,\tau )]$ satisfying
    $$\eqalign{
    V_\alpha^{(-n)}(\pi,\tau)(0) &= 0\, , \cr
    V_\alpha^{(-n)}(\pi,\tau)(q_\alpha (\pi,\tau)) &= v_\alpha (\pi,\tau)\, . \cr}$$
    (In particular $V_\alpha^{(0)}(\pi,\tau)$ is the affine map on $[0,q_\alpha(\pi,\tau )]$ with these boundary values).
\item{$2$.} The vertices of $\Gamma_\alpha^{(-n)}$ are also vertices of $\Gamma_\alpha^{(-n-1)}$.

\noindent
From the behaviour of the coefficients $\Theta_1^{(-n)}$, $\Theta_2^{(-n)}$ and $\Theta^{(-n)}$ it also follows that
\item{$3$.} The sequence $ V_\alpha^{(-n)}(\pi,\tau)$ converges uniformly exponentially fast (with respect to Zorich reparametrized time) at rate $\theta_2$ to a continuous function $V_\alpha (\pi,\tau)$ on
    $[0,q_\alpha(\pi,\tau )]$ (with the same boundary values).
\item{$4$.} The function $V_\alpha (\pi,\tau)$ satisfies a H\"older condition of exponent $\theta$, for any $\theta<\theta_2/\theta_1$.

In property 4, we use the following

\vskip .3 truecm \noindent {\bf Lemma.} {\it For almost all $(\pi,\lambda,\tau)$ we have }
$$
\lim_{n\rightarrow\pm\infty}{1\over Z(n)}\log\hbox{Inf}_\alpha q_\alpha (\pi^{(n)},\tau^{(n)})=0\, .
$$

\vskip .3 truecm\noindent\proof
Let $\gamma$ be a fixed path in $\cD$ such that all coefficients of $B_\gamma$ are positive, and let
$\Theta_\gamma$ be the set of $(\pi,\tau)$ whose itinerary under backwards Rauzy--Veech dynamics ends
with $\gamma$. Then we have
$$
\hbox{Inf}_\alpha q_\alpha (\pi,\tau)\ge C_\gamma>0\, ,
$$
for all $(\pi,\tau)\in\Theta_\gamma$. From the ergodicity of the Rauzy--Veech dynamics (in forward and backward time), for almost all $(\pi,\lambda,\tau)$ there exist monotone sequences
$n_k,n_k'$ going to $+\infty$ and $-\infty$ respectively such that
$$
\lim_{k\rightarrow +\infty}{Z_{n_{k+1}}\over Z_{n_k}}=
\lim_{k\rightarrow +\infty}{Z_{n'_{k+1}}\over Z_{n'_k}}=1\, , \,
(\pi^{(n_k)},\tau^{(n_k)})\in\Theta_\gamma\, , \,
(\pi^{(n'_k)},\tau^{(n'_k)})\in\Theta_\gamma\, .
$$
For $n_k\le n<n_{k+1}$, we have
$$
q(\pi^{(n)},\tau^{(n)})=\Omega^{(n)}B_{\gamma (n_k,n)}q(\pi^{(n_k)},\tau^{(n_k)})\, ,
$$
with
$$
\lim_{k\rightarrow +\infty}{1\over Z(n_k)}\log \Omega^{(n)}=0\, ,
$$
which gives the required estimate.
The case of negative $n$ is similar. \qed

\vskip .3 truecm\noindent
We also define the following function $V_* (\pi,\tau)$: if $\alpha_b$, $\alpha_t$ are the last letter of the bottom, top lines of $\pi$, we set:
$$
V_* (\pi,\tau)(x) = \cases{ V_{\alpha_b} (\pi,\tau)(x) & if $0\le x\le q_{\alpha_b}$, \cr
V_{\alpha_t} (\pi,\tau)(x-q_{\alpha_b})+v_{\alpha_b} & if $q_{\alpha_b}\le x\le q_{\alpha_b}+q_{\alpha_t}$,  \cr}$$
(with $q_{\alpha_b}=q_{\alpha_b}(\pi,\tau)$, etc.).

\vskip .5 truecm \noindent {\bf 3.4.3 The functions $W_\alpha (\pi,\lambda, \tau)$. }
For $\pi,\tau$ as above, $\alpha\in\cA$, $\lambda\in (\R^+)^\cA$, we can perform with respect to the vector $w(\pi,\lambda,\tau) = v(\pi,\tau )-t(\pi,\lambda,\tau)q(\pi,\tau )$ of Section 3.4.1 the same construction that we did for $v(\pi,\tau )$.
We denote by $w^{(-n)}(\pi ,\lambda ,\tau)$ the image of $w(\pi ,\lambda ,\tau)$ under the Zorich cocycle and
we have
$$
w^{(-n)}(\pi ,\lambda ,\tau)= \Theta_2^{(-n)}w(\pi^{(-n)} ,\lambda^{(-n)} ,\tau^{(-n)})\; .
$$
We obtain functions $W_\alpha (\pi,\lambda,\tau)$, $W_* (\pi,\lambda,\tau)$ which are related to the previous ones by
$$\eqalign{
W_\alpha (\pi,\lambda,\tau)(x) &= V_\alpha (\pi,\tau )(x) - t(\pi,\lambda,\tau) x\, , \cr
W_* (\pi,\lambda,\tau)(x) &= V_* (\pi,\tau )(x) - t(\pi,\lambda,\tau) x\, . \cr}$$

\vskip .5 truecm \noindent {\bf 3.4.4 Relation to Birkhoff sums. } Let $\alpha\in\cA$. Denote as above by  $(\beta_0, \beta_1, \ldots )$ the itinerary of $I^{(0)}_\alpha$ with relation to the partition
$I_\beta^{(-n)}$ till its return to $I^{(0)}$.

Consider the Birkhoff sums
$$\eqalign{
S_\alpha q^{(-n)} (i) &= \sum_{j=0}^{i-1} q_{\beta_j}^{(-n)}(\pi,\tau )\, , \cr
S_\alpha w^{(-n)} (i) &= \sum_{j=0}^{i-1} w_{\beta_j}^{(-n)}(\pi, \lambda ,\tau )\, . \cr}$$
We have then by definition of $\Gamma^{(-n)}$ (for $W(\pi, \lambda ,\tau )$)
$$
W_\alpha (S_\alpha q^{(-n)} (i)) = S_\alpha w^{(-n)} (i)\; . $$
If instead we look at the Birkhoff sums
$$\eqalign{
S_\alpha q (i) &= \sum_{j=0}^{i-1} q_{\beta_j}(\pi^{(-n)},\tau^{(-n)} )\, , \cr
S_\alpha w (i) &= \sum_{j=0}^{i-1} w_{\beta_j}(\pi^{(-n)} ,\lambda^{(-n)} ,\tau^{(-n)} )\, , \cr}$$
we will have, in view of the relation between $q^{(-n)},w^{(-n)}$ and $q,w$:
$$\eqalign{
S_\alpha q (i) &= (\Theta_1^{(-n)} )^{-1} S_\alpha q^{(-n)} (i) \, , \cr
S_\alpha w (i) &= (\Theta_2^{(-n)} )^{-1} S_\alpha w^{(-n)} (i)
\, , \cr}$$
hence
$$
S_\alpha w (i) = (\Theta_2^{(-n)} )^{-1} W_\alpha (\Theta_1^{(-n)} S_\alpha q (i)) \, .
$$
In view of this formula one can think of $W_\alpha$ as the ``limit shape'' for the Birkhoff sum of $w$.

\vskip .5 truecm \noindent {\bf 3.4.5 Functional equation. } Here we relate the $W_\alpha (\pi, \lambda ,\tau )$ to the $W_\alpha (\pi^{(-1)},$ $ \lambda^{(-1)}, \tau^{(-1)} )$. The relation is a consequence of the formulas
$$
\eqalign{
q^{(-1)}(\pi,\tau ) &= \Theta_1^{(-1)} q(\pi^{(-1)},\tau^{(-1)} )\, , \cr
w^{(-1)}(\pi,\lambda ,\tau ) &= \Theta_2^{(-1)} w(\pi^{(-1)},\lambda^{(-1)},\tau^{(-1)} )
\, . \cr}
$$
Indeed, {\it if $\alpha$ is not the loser} of the arrow from $\pi^{(-1)}$ to $\pi^{(0)}$, we obtain
$$
W_\alpha (\pi,\lambda ,\tau )(x) = \Theta_2^{(-1)} W_\alpha (\pi^{(-1)},\lambda^{(-1)},\tau^{(-1)} )\left({x\over \Theta_1^{(-1)} }\right) \; .
$$
{\it If $\alpha$ is the loser} of this arrow, we obtain
 $$
W_{\alpha_l} (\pi,\lambda ,\tau )(x) = \Theta_2^{(-1)} W_* (\pi^{(-1)},\lambda^{(-1)},\tau^{(-1)} )\left({x\over \Theta_1^{(-1)} }\right) \; .
$$

\vskip .5 truecm \noindent {\bf 3.5 On the direction of $w$}

Recall that in Section 3.3.1 we want to bound from above the Birkhoff sums of $w$ at some point $x^*$.
In Section 3.4.4 we have related the Birkhoff sums of $w(\pi^{(-n)} ,\lambda^{(-n)} ,\tau^{(-n)} )$
to the limit shape $W_\alpha (\pi, \lambda ,\tau)$. In Section 3.7 the point $x^*$ will be defined using
the maximum of $W_\alpha (\pi^{(n)} ,\lambda^{(n)} ,\tau^{(n)} )$ (for $n>0$ large). Therefore we
need to compare these functions $W_\alpha (\pi^{(n)} ,\lambda^{(n)} ,\tau^{(n)} )$ to their maximum values.
In order to do this, the Proposition below is a crucial technical step.

\vskip .5 truecm \noindent {\bf 3.5.1 } The Rauzy operations $R_t$, $R_b$ in $\cR$ do not change the first
letter of the bottom and top lines of elements of $\cR$. So there is a letter $a\in\cA$ which is the first letter in the top line of any element of $\cal R$. Consider the set $\Upsilon$ of $(\pi,\lambda,\tau)$ with
$\pi\in\cR$, $\lambda\in (\R^+)^\cA$, $\tau\in\Theta_\pi$, which satisfy the following properties
\item{(i)} $a$ is the last letter of the bottom line of $\pi$;
\item{(ii)} $a$ is the loser of the next step of the Rauzy--Veech algorithm for $(\pi,\lambda,\tau)$: if
$\alpha$ is the last letter of the top line of $\pi$, we have $\lambda_\alpha >\lambda_a$;
\item{(iii)} $w_a(\pi,\lambda,\tau)(w_a(\pi,\lambda,\tau)+w_\alpha(\pi,\lambda,\tau))<0$.

Here $w(\pi,\lambda,\tau)$ is the vector associated to the exponent $\theta_2$ defined in 3.4.1. There were two possible choices for $w$ but obviously property (iii) does not depend on this choice. Observe also that there are elements $\pi\in\cR$ satisfying (i): since the Rauzy--Veech expansion of a
standard i.e.m.\ with no connections produces an $\infty$--complete path (see e.g.\ [Y1])
the letter $a$ must be the winner of at least one arrow in $\cal D$ and this can only occur
when $a$ is the last letter of the bottom line.

\vskip .3 truecm\noindent {\bf Proposition.}{\it The set $\Upsilon$ has positive measure.}

\vskip .5 truecm \noindent\proof The rest of this Section 3.5 is devoted to the proof of this assertion.

\vskip .5 truecm \noindent {\bf 3.5.2 } Recall that
$$
w(\pi,\lambda,\tau )= v(\pi,\tau ) - {<\lambda, v> \over <\lambda , q>} q(\pi, \tau )\, .
$$
It will be crucial for the whole argument that the vector w, although it depends only measurably
on $\tau$, is a smooth and explicit function of $\lambda$.

In view of (ii), the vector $\lambda $ is allowed to vary in a convex cone whose extremal vectors $\lambda^{(\beta)}$ are
given by
$$
\eqalign{
\bullet \, & \lambda^{(\beta)}_\gamma := \delta_{\gamma\beta}\, , \;\beta\not= a\, \cr
\bullet \, & \lambda^{(a)}_\gamma := \delta_{\gamma a}+\delta_{\gamma\alpha}\, , \cr}
$$
where $\delta_{\gamma\beta}$, $\delta_{\gamma a}$ and $\delta_{\gamma\alpha}$ denote the Kronecker symbol.
The corresponding values for $w_a$ are
$$
\eqalign{
\bullet \, & v_a - {v_\beta\over q_\beta} q_a \, , \;\beta\not= a\, \cr
\bullet \, & v_a - {v_\alpha +v_a\over q_\alpha +q_a} q_a \, . \cr}
$$
We see that these values have the same sign if and only if ${v_a\over q_a}$ is either larger than all
other ${v_\beta\over q_\beta}$ or smaller than these quantities. Furthermore, if a change of sign of $w_a$
occurs, we want that $w_a+w_\alpha$ does not change sign at the same time, and this occurs if and only if
${v_a\over q_a}={v_\alpha \over q_\alpha}$. We will prove below the following two results

\vskip .3 truecm\noindent {\bf Proposition 1.}{\it Let $\pi\in\cR$ such that $a$ is the first top letter and last bottom letter of
$\pi$. For all $\alpha\in\cA$, $\alpha\not= a$ and almost all $\tau$ we have}
$$
v_a(\pi,\tau )q_\alpha(\pi,\tau )-v_\alpha (\pi,\tau )q_a(\pi,\tau )\not= 0\, .
$$

\vskip .3 truecm\noindent {\bf Proposition 2.}{\it There exist $\pi\in\cR$, with last bottom letter $a$,
letters $b,c$ and a positive measure set of $\tau$ on which }
$$
{v_c\over q_c} < {v_a\over q_a} < {v_b\over q_b}\; .
$$

These two propositions do indeed imply that $\Upsilon$ has positive measure. Let $a,b,c,\pi,\tau$ be as in Proposition 2; almost surely the conclusion of Proposition 1 is also satisfied. We have $w_a(\pi,\lambda,\tau )<0$ if and only if the linear form $l(\lambda )=<\lambda ,v>-{v_a\over q_a}<\lambda ,q>$ is positive and
$w_a(\pi,\lambda,\tau )+ w_\alpha(\pi,\lambda,\tau )<0$ if and only if the linear form $\tilde l(\lambda )=<\lambda ,v>-{v_a+v_\alpha\over q_a+q_\alpha}<\lambda ,q>$ is positive (here $\alpha$ is the last letter in the top line of $\pi$). One has $l(\lambda^{(b)})>0$, $l(\lambda^{(c)})<0$. Moreover, $l$ and $\tilde l$ are
not proportional thus there exists a set of $\lambda$ of positive measure where $l(\lambda )\tilde l(\lambda )<0$. This concludes the proof of the Proposition. \qed

Obviously, {\it the statement obtained from the Proposition in 3.5.1 and the Propositions 1 and 2 in 3.5.2
by exchanging the role of the top and bottom lines are also true}.

\vskip .5 truecm \noindent {\bf 3.5.3 Proof of Proposition 1.} It is based on the {\it twisting property}
of the Rauzy monoid proved by A.\ Avila and M.\ Viana [AV]. Let us recall the content of this property. For $\pi\in\cR$, be the antisymmetric matrix $\Omega (\pi)$ has been defined by
$$
\Omega_{\beta\gamma}(\pi ) = \cases{ 1 & if $\pi_t\beta <\pi_t\gamma\;$, $\;\pi_b\beta>\pi_b\gamma\;$, \cr
-1  & if $\pi_t\beta >\pi_t\gamma\;$, $\;\pi_b\beta <\pi_b\gamma\;$, \cr
0 & otherwise.\cr}
$$
The subspaces $H(\pi )=\hbox{Im}\,\Omega (\pi )$ have dimension $2g$ and are invariant under the Zorich cocycle, which acts symplectically on these subspaces. Let $\pi\in\cR$, $F\subset H(\pi )$ a subspace of dimension $k$, $0<k<2g$, and $F_1^*, \ldots ,F_l^*\subset H(\pi )$ be subspaces of codimension $k$. The twisting property asserts that there exists a loop $\sigma $ of $\cD$ at $\pi$ such that the image of $F$ under the matrix $B_\sigma$ corresponding to $\sigma$ under the Zorich cocycle is transverse to $F_1^* ,\ldots , F_l^*$.

Consider the $2$--dimensional subspace $F(\pi, \tau )$ generated by $q$ and $v$. As it is associated to
the positive Lyapunov exponents $\theta_1>\theta_2$, it is contained in $H(\pi )$ (the Lyapunov exponents on $\R^\cA/H(\pi )$ are equal to zero).

Let $\pi\in\cR$ be such that $a$ is the first top letter and the last bottom letter of $\pi$ and let $\alpha\in\cA$, $\alpha\not= a$. The relation $v_\alpha q_a-v_aq_\alpha =0$ holds if and only if $F(\pi, \tau )$ is {\it not} transverse to the codimension $2$ subspace
$$
\{y\in\R^\cA\, ,\, y_a=y_\alpha=0\,\}\; .
$$
We claim that the intersection $F^*(\alpha )$ of this subspace with $H(\pi )$ is transverse, hence has
codimension $2$ in $H(\pi )$: indeed, let $\nu\in\R^\cA$, $y=\Omega (\pi )\nu$; as $a$ is the first top and
the last bottom letter of $\pi$ we have
$$
y_a = \sum_{\beta\not= a}\nu_\beta\, ,
$$
On the other hand the coefficient of $\nu_a$ in $y_\alpha$ is $-1$. Therefore the linear forms (of the variable $\nu$) $y_a$ and $y_\alpha$ are not proportional and the claim follows.

Therefore, if the conclusion of Proposition 1 for $\pi,\alpha$ does not hold, there exists a set of positive
measure $X\subset\P (\Theta_\pi)$ such that, for $\tau\in X$, the subspace $F(\pi,\tau)$ is not transverse
to $F^*(\alpha)$.

The following Lemma will be proved below.

\vskip .3 truecm \noindent {\bf Lemma.} {\it Let $\pi\in\cR$, $X\subset\P (\Theta_\pi)$ a subset of positive measure. For any $\varepsilon >0$, there exists a loop $\sigma$ of $\cD$ at $\pi$ such that the measure of
$\P (\Theta_\pi)\setminus (^tB_\sigma(X)\cap\P (\Theta_\pi))$ is $<\varepsilon$. }

\vskip .3 truecm\noindent
From the twisting property and the compactness of the Grassmannians, there exist loops $\sigma_1,\ldots ,\sigma_k$ of $\cD$ at $\pi$ such that, for any $2$--dimensional subspace $F_0\subset H(\pi )$, and
any codimension $2$ subspace $F_0^*\subset H(\pi )$,  $F_0^*$ is transverse to at least one of the
$B_{\sigma_i}F_0$.

Let $\varepsilon >0$, and let $\sigma$ be as in the Lemma above. If $\varepsilon >0$ is small enough, there
exists a set of positive measure $Y\subset\P (\Theta_\pi)$ such that, for $\tau\in Y$, $^tB_{\sigma_i}^{-1}
\tau$ belongs to $^tB_\sigma(X)\cap\P (\Theta_\pi)$ for all $1\le i\le k$. Writing $\tau_i =\, {^tB_{\sigma}^{-1}} ^tB_{\sigma_i}^{-1}\tau$, we have $\tau_i\in X$ for $1\le i\le k$; this means that
$F(\pi, \tau_i)$ is not transverse to $F^*(\alpha)$. As the $F$--bundle is invariant under the Rauzy--Veech dynamics, we have that $F(\pi, \tau_i)= B_\sigma B_{\sigma_i}F(\pi,\tau )$; setting $F_0=F(\pi,\tau)$, $F_0^*
=B_{\sigma}^{-1}F^*(\alpha)$, we see that, for $\tau\in Y$, $B_{\sigma_i}F_0$ is not transverse to $F_0^*$
for all $1\le i\le k$, a contradiction. \qed

\vskip .3 truecm\noindent
{\it Proof of the Lemma.}
On $\sqcup_{\tilde\pi\in\cR}\P (\Theta_{\tilde\pi})$, consider the backwards Rauzy--Veech
dynamics $Q_{RV}^{-1}$ and the return map $R_\pi$ to  $\P (\Theta_\pi)$. The branches
of the iterates $R_\pi^n$ are in one-to-one correspondence with the loops of $\cD$
at $\pi$. For $n\ge 0$, $\tau\in \P (\Theta_\pi)$, let $\eta_n(\tau )$ be the domain
of the branch of $R_\pi^n$ which contains $\tau$. For almost all $\tau$, we
have
$$
\cap_{n\ge 0}\eta_n(\tau )=\{\tau\}\, .
$$
It follows that, for any subset of positive measure $X$ of $\P (\Theta_\pi)$, and
almost all $\tau\in X$, we have
$$
\lim_{n\rightarrow +\infty}{\mu (\eta_n(\tau )\cap X)\over\mu (\eta_n(\tau ))}=1\, .
$$
Next, we observe that the distorsion of the Jacobian of $^tB_\sigma^{-1}$ on
$\P (\Theta_\pi)$ (for a loop $\sigma$ of $\cD$ at $\pi$) is controlled by
$$
{\hbox{Max}_{v\in\Theta_\pi\, ,\Vert v\Vert=1} \Vert ^tB_\sigma^{-1}(v)\Vert\over
\hbox{Min}_{v\in\Theta_\pi\, ,\Vert v\Vert=1} \Vert ^tB_\sigma^{-1}(v)\Vert}\, .
$$
Fix a loop $\sigma_0$ of $\cD$ at $\pi$ such that
$$
\overline{^tB_{\sigma_0}^{-1}(\Theta_\pi )}\subset\{0\}\cup\Theta_\pi\, .
$$
If $\sigma$ is a loop at $\pi$ of the form $\sigma=\sigma_0\sigma_1$, then the
distorsion of the Jacobian of $^tB_\sigma^{-1}$ on $\P (\Theta_\pi)$
is bounded by a constant depending only on $\sigma_0$.

On the other hand, as $Q_{RV}^{-1}$ is ergodic, for almost all $\tau\in X$,
there exist infinitely many integers $n$ such that the loop $\sigma$ associated
to $\eta_n(\tau )$ has the form $\sigma=\sigma_0\sigma_1$. Taking such an integer
large enough and applying $^tB_\sigma$, we obtain that the measure of
$\P (\Theta_\pi)\setminus (^tB_\sigma (X)\cap \P (\Theta_\pi))$
can be made arbitrarily small. \qed

\vskip .3 truecm\noindent {\bf Remark.} Proposition 1 is in general false if we replace $a,\alpha $ by
any two distinct letters: consider in genus $2$
$$
\pi = \left( \matrix{ A & B & C & D & E\cr D & E & C & B & A\cr}\right)\, .
$$
Obviously we have $\{ u_D=u_E\}$ as equation of $H(\pi )$, hence $q_Dv_E-q_Ev_D\equiv 0$.

\vskip .5 truecm \noindent {\bf 3.5.4 Proof of Proposition 2.} Let $c$ be the first letter of the bottom line of all elements of $\cR$: we have $c\not= a$; let $b$ be any letter distinct from $a$ and $c$.
We will prove the inequalities of Proposition 2 up to exchanging $b$ and $c$ (which leaves invariant the
statement of Proposition 2). Let $\pi_0\in\cR$ such
that the last top and bottom letters are $c,a$ respectively (if $\pi\in\cR$ is such that $a$ is the last letter of the bottom line such a $\pi_0$ is obtained by a suitable number of iterations of the Rauzy operation $R_b$). Consider in $\cD$ the subdiagram obtained by erasing the arrows whose winner is not $a,b$ or $c$ and then keeping the connected component $\cD '$ of $\pi_0$. It is easily seen to have the typical form shown in the figure (see [AV], [AGY])

\vskip .5 truecm
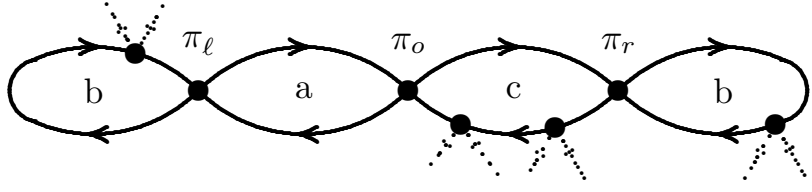
\centerline{\hbox{\psfig{figure=newfigure.epsi,height=2cm}}}

\vskip .5 truecm\noindent
(i.e.\ it is essentially the Rauzy diagram with $d=3$, (see e.g.\ [Y1]) with some meaningless vertices added; only $\pi_0$,
$\pi_l$ and $\pi_r$ have two arrows going out).

For paths contained in $\cD '$, the $a,b,c$ coordinates of vectors are changed under the Zorich cocycle
exactly as in the Rauzy diagram with $d=3$. Consider the vectors in the right halfplane:
$$
\eqalign{
u_a & = u_a (\pi_0,\tau) = (q_a(\pi_0 , \tau ), v_a (\pi_0 ,\tau ))\, ,\cr
u_b & = u_b (\pi_0,\tau) = (q_b(\pi_0 , \tau ), v_b (\pi_0 ,\tau ))\, ,\cr
u_c & = u_c (\pi_0,\tau) = (q_c(\pi_0 , \tau ), v_c (\pi_0 ,\tau ))\,  .\cr}
$$
By Proposition 1 (and its symmetric statement obtained by exchanging top and bottom), for almost all $\tau$, no two among these $3$ vectors are collinear (indeed, $c$ has the same properties than $a$).

If there is a set of $\tau$ of positive measure such that $u_a$ is between $u_b$ and $u_c$ in the right halfplane, the conclusion of Proposition 2 is satisfied; assume therefore that it is not the case.

Next assume that on a set of positive measure the vector $u_a+u_c$ is between $u_a$ and $u_b$. Consider the path $\sigma$ starting at $\pi_0$, going to $\pi_l$ and making $N$--times the $b$--loop at $\pi_l$; the effect on the vectors is the following (we have for each arrow to add the winning vector to the losing one):
$$
\eqalign{
u_a & \longrightarrow u_a' = u_a+Nu_b \, , \cr
u_b & \longrightarrow u_b' = u_b \, , \cr
u_c & \longrightarrow u_c' = u_a+u_c \, . \cr}
$$
If $N$ is large enough then $u_a'$ is between $u_b'$ and $u_c'$ hence the conclusion of Proposition 2 is again satisfied (at $\pi_l$).

Finally, in the remaining case, we would have that, for almost all $\tau$, $u_b$ is between $u_a$ and $u_a+u_c$; the loop at $\pi_0$ obtained by going to $\pi_r$, making $N$ times the $b$--loop at $\pi_r$ and coming back to $\pi_0$ has for effect:
$$
\eqalign{
u_a & \longrightarrow u_a'' = u_a+u_c \, , \cr
u_b & \longrightarrow u_b'' = u_c+(N+1)u_b \, , \cr
u_c & \longrightarrow u_c'' = u_c+Nu_b \, . \cr}
$$
For large $N$, $u_c''$ is between $u_a''$ and $u_b''$, which contradicts the assumption. The proof of Proposition 2 is now complete. \qed

\vskip .5 truecm \noindent {\bf 3.6 Consequences for limit shapes}

\vskip .5 truecm \noindent {\bf 3.6.1 } Let $(\pi ,\lambda ,\tau )$ be a typical point for the Rauzy--Veech dynamics, let $\alpha\in\cA$, and let $W_\alpha (\pi,\lambda ,\tau )$ be the limit shape defined in Section 3.4.3.

\vskip .3 truecm\noindent {\bf Proposition }{\it The extremal values of $W_\alpha (\pi, \lambda ,\tau )$
(minimum and maximum) are not taken at the endpoints of the interval of definition $[0,q_\alpha (\pi ,\tau )]$ of $W_\alpha (\pi, \lambda ,\tau )$.}

\vskip .3 truecm\noindent\proof As the set $\Upsilon$ of the proposition in 3.5.1 has positive measure and the invariant measure for Rauzy--Veech dynamics is conservative and ergodic, there exists (for almost all $(\pi, \lambda ,\tau )$) a positive integer $N$ such that $(\pi^{(-N)}, \lambda^{(-N)} ,\tau^{(-N)} )$ belongs to $\Upsilon$ and the interval $I^{(0)}$ is contained in the first subinterval $I_a^{(-N+1)}$ of $I^{(-N+1)}$. We have then
$$
\eqalign{
W_\alpha (\pi, \lambda ,\tau) (q_a^{(-N)} (\pi, \tau )) &= w_a^{(-N)}(\pi, \lambda ,\tau)\, , \cr
W_\alpha (\pi, \lambda ,\tau) (q_a^{(-N+1)} (\pi, \tau )) &= w_a^{(-N+1)}(\pi, \lambda ,\tau)\, , \cr}
$$
with
$$
\eqalign{
q_a^{(-N+1)} (\pi, \tau ) &= q_a^{(-N)} (\pi, \tau )+q_\alpha^{(-N)}
(\pi, \tau )\, , \cr
w_a^{(-N+1)}(\pi, \lambda ,\tau) &= w_a^{(-N)}(\pi, \lambda ,\tau)+ w_\alpha^{(-N)}(\pi, \lambda ,\tau)\, , \cr}
$$
$\alpha$ being the winner of the arrow from $\pi^{(-N)}$ to $\pi^{(-N+1)}$. By the definition of $\Upsilon$ we have that
$$
w_a^{(-N)}(\pi, \lambda ,\tau) w_a^{(-N+1)}(\pi, \lambda ,\tau) <0
$$
and therefore $0$ is not an extremal value of $W_\alpha (\pi, \lambda ,\tau)$. The other endpoint is treated in a similar manner, exchanging the top and the bottom lines. \qed

\vskip .5 truecm \noindent {\bf 3.6.2 Smallest concave majorant. } Let $F\, : [a,b]\rightarrow \R$ be continuous. The infimum of concave majorants of $F$ on $[a,b]$ is the smallest concave majorant of $F$ and will be denoted by $\hat F$; it is continuous and satisfies $\hat F(a)=F(a)$, $\hat F(b)=F(b)$; moreover,
the maximum values of $F$ and $\hat F$ are the same. We write
$\hat F_r'$, $\hat F_l'$ for the right and left derivatives of $\hat F$.

\vskip .3 truecm\noindent {\bf Proposition }{\it Let $(\pi, \lambda ,\tau )$ be a typical point for Rauzy--Veech dynamics and let $\alpha\in\cA$. We have}
$$
\eqalign{
\hat W_{\alpha ,r}' (\pi, \lambda ,\tau ) (0) &= +\infty\, , \cr
\hat W_{\alpha ,l}' (\pi, \lambda ,\tau ) (q_\alpha (\pi, \tau )) &= -\infty\, , \cr
\hat W_{* ,r}' (\pi, \lambda ,\tau ) (q_{\alpha_b} ) &= \hat W_{* ,l}' (\pi, \lambda ,\tau ) (q_{\alpha_b} )\not= 0
\, . \cr}
$$

\vskip .3 truecm\noindent\proof The first two assertions are a very slight extension of the Proposition in 3.6.1: in the proof of this proposition we first replace the set $\Upsilon$ of Section 3.5.1 by the slightly smaller set $\Upsilon_\delta $ obtained by replacing condition (iii) in 3.5.1 by
$$
 w_a(w_a+w_\alpha )<0\, , \;\hbox{and}\, |w_a|>\delta \, \hbox{and}\, |w_a+w_\alpha |>\delta\; . \leqno\hbox{(iii)}_\delta
$$
If $\delta >0$ is small enough, this has still positive measure. Now, the integer $N$ in the proof of Proposition 3.6.1 can be taken arbitrarily large; as $q_a^{(-N)}$ and $w_a^{(-N)}$ go down exponentially fast (in Zorich time) at respective rates $\theta_1>\theta_2$, this implies the first two assertions of the Proposition.

For the last assertion, it follows from the definition of $V_*$ and the first two assertions that we have
$$
\hat V_*(\pi, \tau )(q_{\alpha_b})>V_*(\pi, \tau ) (q_{\alpha_b})\, .
$$
It follows that $\hat V_*$ is affine in a neighborhood of $q_{\alpha_b}$, in particular
$\hat V_{*,r}'(q_{\alpha_b})=\hat V_{*,l}'(q_{\alpha_b})$.

Now, obviously we have
$$
\hat W_* (\pi, \lambda ,\tau )(x) = \hat V_*(\pi ,\tau )(x)-{<\lambda , v>\over <\lambda,q>}x\, ,
$$
(adding an affine function to $F$ adds the same affine function to the smallest concave majorant).
Therefore we have
$$
\hat W_*' (\pi, \lambda ,\tau )(q_{\alpha_b})=0
$$
if and only if
$$
{<\lambda , v>\over <\lambda,q>} = \hat V_*'(\pi ,\tau )(q_{\alpha_b})
$$
which has $\lambda$--measure zero for any given $(\pi ,\tau )$. \qed

\vskip .5 truecm \noindent {\bf 3.6.3 Corollary. } {\it The function $W_\alpha (\pi, \lambda ,\tau )$
takes its maximum value at a unique point $x_\alpha^{\hbox{\sevenrm max}}(\pi, \lambda ,\tau )$ (for almost all
$(\pi, \lambda ,\tau )$). }

\vskip .5 truecm\noindent\proof
Let $(\pi, \lambda ,\tau )$ be a typical point and $\alpha\in\cA$. By the functional equation of Section 3.4.5, $W_\alpha (\pi, \lambda ,\tau )$ is a rescaled version of either $W_\alpha (\pi^{(-1)}, \lambda^{(-1)} ,\tau^{(-1)} )$ (if $\alpha$ is not the loser of the arrow from $\pi^{(-1)}$ to $\pi$) or $W_* (\pi^{(-1)}, \lambda^{(-1)} ,\tau^{(-1)} )$ (if $\alpha$ is the loser of this arrow).

In this last case, by the last assertion of Proposition 3.6.2, $W_* (\pi^{(-1)},$ $\lambda^{(-1)} ,\tau^{(-1)} )$ does not take its maximum value both in $[0,q_{\alpha_b}(\pi^{(-1)},\tau^{(-1)} )]$ and in $[q_{\alpha_b}(\pi^{(-1)},\tau^{(-1)} ), q_{\alpha_b}(\pi^{(-1)},\tau^{(-1)} )+q_{\alpha_t}(\pi^{(-1)},\tau^{(-1)} )]$ (otherwise we would have $\hat W_*'(\pi^{(-1)}, \lambda^{(-1)} ,\tau^{(-1)} )(q_{\alpha_b}(\pi^{(-1)},\tau^{(-1)} ))=0$).

In view of the definition of $W_*$, this means that the set ${\cal M}$ where $W_\alpha (\pi, \lambda ,\tau )$
takes its maximum value is a rescaled version of the set where $W_{\alpha (1)} (\pi^{(-1)}, \lambda^{(-1)} ,$ $\tau^{(-1)} )$ takes its maximum value, for some $\alpha (1)\in\cA$. Iterating this procedure, we obtain that ${\cal M}$ is a rescaled version (by a factor $\Theta_1^{(-n)}$) of the set where $W_{\alpha (n)}(\pi^{(-n)}, \lambda^{(-n)} ,\tau^{(-n)} )$ takes its maximum value, for some letter $\alpha (n)\in\cA$. As the $q_\alpha$
are bounded by $1$ this proves that the diameter of ${\cal M}$ is smaller than $\Theta_1^{(-n)}$ for all
$n\ge 0$, hence it is a point.  \qed

\vskip .3 truecm\noindent
A similar result is true for minimum values. The function $W_* (\pi, \lambda ,\tau )$ also takes its maximum value at a unique point $x_*^{\hbox{\sevenrm max}} (\pi, \lambda ,\tau )$. By the proposition in 3.6.1 we know that
$x_*^{\hbox{\sevenrm max}} (\pi, \lambda ,\tau )$ is distinct from $0$, $q_{\alpha_b}$, $q_{\alpha_b}+q_{\alpha_t}$.
Observe that we have
$$
\eqalign{
x_*^{\hbox{\sevenrm max}} (\pi, \lambda ,\tau ) \in (0,q_{\alpha_b}) &\Longleftrightarrow \hat W_*' (\pi, \lambda ,\tau )(q_{\alpha_b}) <0\, , \cr
x_*^{\hbox{\sevenrm max}} (\pi, \lambda ,\tau ) \in (q_{\alpha_b}, q_{\alpha_b}+q_{\alpha_t})
&\Longleftrightarrow \hat W_*' (\pi, \lambda ,\tau )(q_{\alpha_b}) >0\, . \cr}
$$
Assume for instance that $x_*^{\hbox{\sevenrm max}} (\pi, \lambda ,\tau ) \in (0,q_{\alpha_b})$. As $W_*$ and $\hat W_*$ coincide at $x_*^{\hbox{\sevenrm max}}$, we have, for $x\in [q_{\alpha_b}, q_{\alpha_b}+q_{\alpha_t}]$
$$
\eqalign{
W_*(x) &\le \hat W_*(x) \cr &\le \hat W_* (q_{\alpha_b}) + \hat W_*' (q_{\alpha_b})
(x-q_{\alpha_b})\cr &\le W_* (x_*^{\hbox{\sevenrm max}})+ \hat W_*' (q_{\alpha_b})
(x-q_{\alpha_b}) \; . \cr}
$$
This will provide a satisfactory control of $W_*$ if $|\hat W_*' (q_{\alpha_b})|$ is not too small and
$(x-q_{\alpha_b})$ is not too small. When $x$ is very close to $q_{\alpha_b}$, we will rely on a direct control on $W_* (x_*^{\hbox{\sevenrm max}})-W_*(q_{\alpha_b})$, based on the Proposition in 3.5.1.

\vskip .5 truecm \noindent {\bf 3.7 Proof of the Proposition in 3.3.1}

\vskip .5 truecm \noindent {\bf 3.7.1 } Let $(\pi, \lambda ,\tau )$ be a typical point for the Rauzy--Veech dynamics.

We observe first that, if $\tilde w$ is a vector in the subspace $E_2(\pi ,\lambda )$, Zorich has proved [Z2] that the Birkhoff sums $S_n\tilde w$ satisfy, uniformly on $I^{(0)}$, an estimate
$$
\Vert S_n (\tilde w)\Vert_{C^0}\le C(\varepsilon )|n|^{\omega +\varepsilon }\, ,
$$
for all $\varepsilon >0$; here $\omega$ is either $0$ if $g=2$ or $\theta_3/\theta_1$ if $g\ge 3$. In any case, we have $\omega+\varepsilon < \theta_2/\theta_1 -\varepsilon$ for small $\varepsilon$, hence the order is smaller than the one in Proposition 3.3.1.

It follows that it is sufficient to prove the estimate of Proposition 3.3.1 when $w$ is
``the'' vector $w(\pi, \lambda ,\tau )$ considered above (there are actually two vectors to consider, opposite to each other).

\vskip .5 truecm \noindent {\bf 3.7.2 } Recall the relation between Birkhoff sums and limit shapes from
Section 3.4.4:
$$
S_\alpha w(i) = \Theta_2^{(n)} W_\alpha ((\Theta_1^{(n)})^{-1}S_\alpha q (i))\, ,
$$
where
\item{$\bullet$} $S_\alpha q(i)=\sum_{j=0}^{i-1} q_{\beta_j}(\pi ,\tau )\, $,
\item{$\bullet$} $ S_\alpha w(i) =\sum_{j=0}^{i-1} w_{\beta_j}(\pi ,\lambda ,\tau )\, $,
\item{$\bullet$} $\beta_0,\beta_1,\ldots$ is the itinerary of $I_\alpha^{(n)}$ with relation to the partition $I_\beta^{(0)}$ of $I^{(0)}$,
\item{$\bullet$} $W_\alpha=W_\alpha (\pi^{(n)}, \lambda^{(n)}, \tau^{(n)})$ is the limit shape at $(\pi^{(n)}, \lambda^{(n)}, \tau^{(n)})$,
\item{$\bullet$} the real number $\Theta_1^{(n)}=\Theta_1^{(n)}(\pi,\lambda,\tau)>0$ is defined by the relation $q^{(n)}(\pi , \tau)=\Theta^{(n)}_1 q
(\pi^{(n)}, \tau^{(n)} )$ where $q^{(n)}(\pi ,\tau )$ is the image of $q(\pi ,\tau )$ under the Zorich cocycle,
\item{$\bullet$} the real number $\Theta_2^{(n)}=\Theta_2^{(n)}(\pi,\lambda,\tau)$ is similarly defined by $w^{(n)}(\pi ,\lambda ,\tau )=\Theta_2^{(n)}
w (\pi^{(n)}, \lambda^{(n)}, \tau^{(n)} )$,
\item{$\bullet$} $i$ varies from $0$ to the return time of $I_\alpha^{(n)}$ in $I^{(n)}$ under $T^{(0)}$.

\noindent
We assume that the choices of signs for $w(\pi,\lambda, \tau )$ and $w( \pi^{(n)}, \lambda^{(n)}, \tau^{(n)} )$
are such that
$$
\Theta_2^{(n)} >0\; .
$$
By Corollary 3.6.3, for almost all $(\pi, \lambda ,\tau )$, all $\alpha\in\cA$, all $n\ge 0$,
$W_\alpha (\pi^{(n)}, \lambda^{(n)}, \tau^{(n)} )$ has a unique maximum at some $x_\alpha^{\hbox{\sevenrm max}}=
x_\alpha^{\hbox{\sevenrm max}}(\pi^{(n)}$ , $\lambda^{(n)}$ , $\tau^{(n)} )$. Let $i$ be the integer such that
$$
S_\alpha q(i) < \Theta_1^{(n)}x_\alpha^{\hbox{\sevenrm max}} <S_\alpha q(i+1)\, , \leqno(4)
$$
where the inequalities are strict, by Proposition 3.6.1.

Let $I_\alpha^{\hbox{\sevenrm max}}(n)$ be the image of $I_\alpha^{(n)}$ by $(T^{(0)})^i$.

Consider what happens when going from  $n$ to $n+1$. If $\alpha$ is {\it not the loser} of the arrow from $\pi^{(n)}$ to $\pi^{(n+1)}$, $W_\alpha (\pi^{(n+1)}, \lambda^{(n+1)}, \tau^{(n+1)} )$ is a rescaled version of $W_\alpha (\pi^{(n)}, \lambda^{(n)}, \tau^{(n)} )$, hence the respective maxima correspond. Therefore the values of $i$ are the same, and $I_\alpha^{\hbox{\sevenrm max}}(n+1)$ is equal to (if $\alpha$ is not the winner)
or contained in (if $\alpha $ is the winner) $I_\alpha^{\hbox{\sevenrm max}}(n)$ (because $I_\alpha^{(n+1)}$ is equal to, resp.\ contained in, $I_\alpha^{(n)}$).

If $\alpha$ {\it is the loser} of the arrow from $\pi^{(n)}$ to $\pi^{(n+1)}$, $W_\alpha (\pi^{(n+1)}, \lambda^{(n+1)}, \tau^{(n+1)} )$ is a rescaled version of $W_* (\pi^{(n)}, \lambda^{(n)}, \tau^{(n)} )$. Write as usual $\alpha_b$ (resp.\ $\alpha_t$) for the last letters in the bottom (resp.\ top) lines of $\pi^{(n)}$. The maximum $x_*^{\hbox{\sevenrm max}}$ is either $x_{\alpha_b}^{\hbox{\sevenrm max}}$ or $q_{\alpha_b}+x_{\alpha_t}^{\hbox{\sevenrm max}}$; in the first case, the values of $i$ for $I_\alpha^{\hbox{\sevenrm max}}(n+1)$ and $I_{\alpha_b}^{\hbox{\sevenrm max}}(n)$ are again the same, and $I_\alpha^{(n+1)}$ is a subinterval of $I_{\alpha_b}^{(n)}$, hence $I_\alpha^{\hbox{\sevenrm max}}(n+1)\subset I_{\alpha_b}^{\hbox{\sevenrm max}}(n)$; in the second case, the values of $i$ for $I_\alpha^{\hbox{\sevenrm max}}(n+1)$ and $I_{\alpha_t}^{\hbox{\sevenrm max}}(n)$  differ by the return time $Q$ of $I_{\alpha_b}^{(n)}$ in $I^{(n)}$, and the image of $I_{\alpha}^{(n+1)}$ under $(T^{(0)})^Q$ is contained in $I_{\alpha_t}^{(n)}$, hence $I_\alpha^{\hbox{\sevenrm max}}(n+1)$ is contained in $I_{\alpha_t}^{\hbox{\sevenrm max}}(n)$.

Thus, we have the following

\vskip .3 truecm \noindent {\bf Lemma. } {\it For each $n$, the intervals $I_\alpha^{\hbox{\sevenrm max}}(n)$
are disjoint. They satisfy
$$
I_\alpha^{\hbox{\sevenrm max}}(n+1)\subset I_{\eta_n(\alpha )}^{\hbox{\sevenrm max}}(n)
$$
where $\eta_n(\alpha )=\alpha$ except possibly when $\alpha$ is the loser of the arrow from $\pi^{(n)}$ to $\pi^{(n+1)}$; in this case $\eta_n(\alpha )$ is either $\alpha$ or the winner of the same arrow. }

\vskip .3 truecm\noindent\proof The last assertion has been proved above, the first one is clear because the orbits of the $I_{\alpha}^{(n)}$ are disjoint till their return time. \qed

We can now specify the point $x^*$ in Proposition 3.3.1. Indeed, take any sequence $(\alpha_n)_{n\ge 0}\subset\cA$ such that
$$
\eta_n (\alpha_{n+1})=\alpha_n\, .
$$

\vskip .3 truecm \noindent {\bf Remark. } It is reasonable to expect that for almost all $(\pi, \lambda ,\tau )$ such a sequence is unique.

\vskip .3 truecm \noindent
The point $x^*$ is defined to be
$$
x^* = \cap_{n\ge 0} \overline{I_{\alpha_n}^{\hbox{\sevenrm max}}(n)}\; .
$$

\vskip .5 truecm \noindent {\bf 3.7.3 } The Birkhoff sums of $w$ at $x^*$ and the functions $W_\alpha$ are related as follows.

Denote by $Q^+(n)\ge 0$ (respectively $Q^-(n)\le 0$) the first entrance time in the future (resp.\ in the past)
of $x^*$ in $I^{(n)}$ under $T^{(0)}$. The sequence $Q^{+}(n)$ is non decreasing and the sequence
$Q^-(n)$ is non increasing.

Moreover, for almost all $(\pi, \lambda ,\tau )$, one has $(\pi^{(n)}, \lambda^{(n)}, \tau^{(n)} )\in\Upsilon_\delta$ for infinitely many $n\ge 0$, where $\Upsilon_\delta$ is the set defined in 3.6.2. It follows that there are arbitrarily large values of $n$ such that the maximum $x^{\hbox{\sevenrm max}}_{\alpha_n}(\pi^{(n)}, \lambda^{(n)}, \tau^{(n)} )$ of
$W_{\alpha_n}(\pi^{(n)}, \lambda^{(n)}, \tau^{(n)} )$ is not exponentially small w.r.t.\ Zorich time $Z(n)$. This implies that the integer $i$ in formula (4) above goes to $+\infty$ and thus
$$
\lim_{n\rightarrow +\infty}Q^-(n) = -\infty \, ,
$$
and similarly one has
$$
\lim_{n\rightarrow +\infty}Q^+(n) = +\infty \, .
$$

Given some integer $j$, we want to estimate the Birkhoff sum $S_jw(x^*)$.

Assume for instance that $j$ is positive (the other case is symmetric) and let $n$ be such that
$$
Q^+(n)<j\le Q^+(n+1)\, .
$$
For $m\ge 0$, let $i_m\ge 0$ be the integer such that
$$
I_{\alpha_m}^{\hbox{\sevenrm max}}(m)=T^{i_m} (I_{\alpha_m}^{(m)})\, .
$$
\vskip .3 truecm\noindent
{\bf Claim.}{ \it $\alpha_{n+1}$ is the loser of the arrow from $\pi^{(n)}$ to $\pi^{(n+1)}$.}

\vskip .3 truecm\noindent\proof
Assume that this is not the case. Then the discussion before the lemma in Section 3.7.2 shows that
$i_n=i_{n+1}$; on the other hand, the return times of $I_{\alpha_{n+1}}^{(n)}$ in $I^{(n)}$ and
$I_{\alpha_{n+1}}^{(n+1)}$ in $I^{(n+1)}$ are the same. Then we would have $Q^+(n)=Q^+(n+1)$, a contradiction. \qed

\vskip .3 truecm\noindent
It follows from the claim that $W_{\alpha_{n+1}}(\pi^{(n+1)}, \lambda^{(n+1)}, \tau^{(n+1)} )$ is a rescaled
version of $W_*=W_*(\pi^{(n)}, \lambda^{(n)}, \tau^{(n)} )$.

We have then
$$
\eqalign{
S_jw(x^*) &= \sum_{k=i_{n+1}}^{i_{n+1}+j-1}w_{\beta_k}(\pi ,\lambda ,\tau )\, \cr
&= \Theta_2^{(n)} \left(W_*([\Theta_1^{(n)}]^{-1}S_{\alpha_{n+1}}q(i_{n+1}+j)) -
W_*([\Theta_1^{(n)}]^{-1}S_{\alpha_{n+1}}q(i_{n+1}))\right)\, . \cr}
$$

\vskip .3 truecm\noindent
{\bf Claim.}{ \it We have $i_{n+1}=i_n$, $\pi_b^{(n)}(\alpha_{n+1})=d$ and
$$
[\Theta_1^{(n)}]^{-1}S_{\alpha_{n+1}}q(i_{n+1})\in [0,q_{\alpha_b}(\pi^{(n)}, \lambda^{(n)}, \tau^{(n)} )]
$$
and }
$$
\eqalign{
[\Theta_1^{(n)}]^{-1}S_{\alpha_{n+1}}q(i_{n+1}+j) &\in [q_{\alpha_b}(\pi^{(n)}, \lambda^{(n)}, \tau^{(n)} ),\cr
& q_{\alpha_b}(\pi^{(n)}, \lambda^{(n)}, \tau^{(n)} )+q_{\alpha_t}(\pi^{(n)}, \lambda^{(n)}, \tau^{(n)} )]\, .\cr}
$$

\vskip .3 truecm\noindent \proof
We refer again to the discussion before the lemma in Section 3.7.2. We claim that in this discussion we must have that $x_*^{\hbox{\sevenrm max}}$ is $x_{\alpha_b}^{\hbox{\sevenrm max}}$. (Otherwise this discussion
shows that $Q^+(n+1)=Q^+(n)$). We have seen in Section 3.7.2 that then we have $i_n=i_{n+1}$, $\pi_b^{(n)}(\alpha_{n+1})=d$ and thus
$$
[\Theta_1^{(n)}]^{-1}S_{\alpha_{n+1}}q(i_{n+1})\in [0,q_{\alpha_b}(\pi^{(n)}, \lambda^{(n)}, \tau^{(n)} )]\, .
$$
Moreover, we have
$$
[\Theta_1^{(n)}]^{-1}S_{\alpha_{n+1}}q(i_{n+1}+Q^+(n))=q_{\alpha_b}(\pi^{(n)}, \lambda^{(n)}, \tau^{(n)} )\, , $$
and
$$
[\Theta_1^{(n)}]^{-1}S_{\alpha_{n+1}}q(i_{n+1}+Q^+(n+1))=q_{\alpha_b}(\pi^{(n)}, \lambda^{(n)}, \tau^{(n)} )+q_{\alpha_t}(\pi^{(n)}, \lambda^{(n)}, \tau^{(n)} )\, , $$
Hence
$$
\eqalign{
[\Theta_1^{(n)}]^{-1}S_{\alpha_{n+1}}q(i_{n+1}+j)&\in [q_{\alpha_b}(\pi^{(n)}, \lambda^{(n)}, \tau^{(n)} ),\cr &q_{\alpha_b}(\pi^{(n)}, \lambda^{(n)}, \tau^{(n)} )+q_{\alpha_t}(\pi^{(n)}, \lambda^{(n)}, \tau^{(n)} )]\, .
\cr}
$$
This concludes the proof of the claim. \qed

\vskip .3 truecm\noindent

Let
$$
\eqalign{
y^\dagger &= (\Theta_1^{(n)})^{-1}S_{\alpha_{n+1}}q(i_{n+1}+j)\, , \cr
y^* &= (\Theta_1^{(n)})^{-1}S_{\alpha_{n+1}}q(i_{n+1})\, . \cr}
$$
From the construction of $W_*$ we have
$$
|\Theta_2^{(n)}(W_*(y^*)-W_*(x_*^{\hbox{\sevenrm max}}))|\le C\, ,
$$
where the majorant $C$ depends on $(\pi , \lambda ,\tau )$ but {\it not} on $n$. We therefore are left with the estimation of
$$
\Theta_2^{(n)} (W_*(y^\dagger)-W_*(x_*^{\hbox{\sevenrm max}}))\, ,
$$
when  $x_*^{\hbox{\sevenrm max}}\in [0,q_{\alpha_b}]$, $y^\dagger\in [q_{\alpha_b},
q_{\alpha_b}+q_{\alpha_t}]$.

\vskip .5 truecm \noindent {\bf 3.7.4 }
For $n\ge 0$, write $W_*^{\hbox{\sevenrm max}}(n)$ for the maximum value of $W_*(\pi^{(n)}, \lambda^{(n)}, \tau^{(n)} )$ in its domain $[0, q_{\alpha_b}+q_{\alpha_t}]$. If the maximum value is taken in $[0,q_{\alpha_b}]$, let $\tilde W_*^{\hbox{\sevenrm max}}(n)$ be the maximum value of $W_*$ in
$[q_{\alpha_b}, q_{\alpha_b}+q_{\alpha_t}]$; if the maximum value of $W_*$ is taken in $[q_{\alpha_b},
q_{\alpha_b}+q_{\alpha_t}]$ , let $\tilde W_*^{\hbox{\sevenrm max}}(n)$ be the maximum value in $[0, q_{\alpha_b}]$.

To complete the proof of Proposition 3.3.1, it is therefore sufficient to prove the following estimate:

\vskip .3 truecm \noindent {\bf Proposition. } {\it For almost all $(\pi, \lambda ,\tau )$ one has
$$
\lim_{n\rightarrow +\infty} {1\over Z(n)}\log (W_*^{\hbox{\sevenrm max}}(n)-\tilde W_*^{\hbox{\sevenrm max}}(n))=0\, ,
$$
where $Z(n)$ is the  Zorich time defined in Section 1.4.}

\vskip .3 truecm\noindent\proof
Let us first deal with the upper bound (which is actually not needed for our purposes). From the normalisation of $v(\pi ,\tau )$ and the construction of $V_*(\pi ,\tau )$, it is clear that, for almost all
$\pi, \lambda, \tau)$, we have
$$
\limsup_{n\rightarrow +\infty} {1\over Z(n)}\log\Vert V_*(\pi^{(n)}, \lambda^{(n)},\tau^{(n)})\Vert_{C^0}
\le 0\; .
$$
On the other hand, from the Lemma in 3.4.2, we also have almost surely
$$
\limsup_{n\rightarrow +\infty} {1\over Z(n)}\log |t (\pi^{(n)}, \lambda^{(n)},\tau^{(n)})|\le 0\; .
$$
It follows that
$$
\limsup_{n\rightarrow +\infty} {1\over Z(n)}\log\Vert W_*(\pi^{(n)}, \lambda^{(n)},\tau^{(n)})\Vert_{C^0}
\le 0\; .
$$
For the lower bound, we will deal first with a neighborhood of the central point $q_{\alpha_b}$ (of
a size which is not exponentially small in $Z(n)$).

We apply Birkhoff ergodic theorem to the Rauzy--Veech dynamics (in Zorich time) and to the characteristic function of the set $\Upsilon_\delta $. We see that for any $n$ there exists $n'<n$ such that $I^{(n)}$ is contained in the first interval $I_a^{(n'+1)}$, $(\pi^{(n')}, \lambda^{(n')}, \tau^{(n')} )$ belongs to $\Upsilon_\delta $, and the ratio
${Z(n)-Z(n')\over Z(n)}$ converges to $0$ as $n\rightarrow +\infty$.
By definition of $\Upsilon_\delta$ and the scaling rules, there exists a point $x_1\in [q_{\alpha_b},
q_{\alpha_b}+q_{\alpha_t}]$ such that
$$
W_*(x_1)-W_*(q_{\alpha_b})\ge \omega (n):= \delta {\hbox{Min}\, [\Theta_2^{(n')},\Theta_2^{(n'+1)}]\over \Theta_2^{(n)} }\; .
$$
From the properties of $n'$ and $\Theta_2^{(n')},\Theta_2^{(n)} $ we have
$$
\lim_{n\rightarrow +\infty}{1\over Z(n)}\log\omega (n)=0\; .
$$
We therefore have
$$
W_*^{\hbox{\sevenrm max}}(n)-W_*(q_{\alpha_b})\ge\omega(n)\; .
$$
We have seen in the Lemma in 3.4.2 that
$$
\lim_{n\rightarrow +\infty}{1\over Z(n)}\log\hbox{Min}_\alpha q_\alpha (\pi^{(n)}, \tau^{(n)})=0\, .
$$
It follows that
$$
|W_*(y)-W_*(q_{\alpha_b})|\le{1\over 2}\omega (n)
$$
for $|y-q_{\alpha_b}|<r(n)$, where $r(n)$ can be chosen so that
$$
\lim_{n\rightarrow +\infty}{1\over Z(n)}\log r(n)=0\, .
$$
We have obtained so far that
$$
W_*^{\hbox{\sevenrm max}}(n)-W_*(y)\ge {1\over 2}\omega(n)\; ,
$$
for $|y-q_{\alpha_b}|<r(n)$. We will deal now with the case where $|y-q_{\alpha_b}|>r(n)$.
In this case, as $y$ and $x_*^{\hbox{\sevenrm max}}$ are separated by $q_{\alpha_b}$, we use the
smallest concave majorant  $\hat W_*$ of Section 3.6.2 to get
$$
\eqalign{
W_*^{\hbox{\sevenrm max}}(n)-W_*(y) &=W_*(x_*^{\hbox{\sevenrm max}})-W_*(y)=
\hat W_*(x_*^{\hbox{\sevenrm max}})-W_*(y)\cr
&\ge \hat W_*(q_{\alpha_b})-\hat W_*(y) \ge |\hat W_*'(q_{\alpha_b})|r(n)\, . \cr}
$$
To complete the proof of the Proposition we use the

\vskip .3 truecm\noindent {\bf Claim} {\it For almost all $(\pi , \lambda , \tau )$ we have}
$$
\liminf_{n\rightarrow +\infty}{1\over n}\log |\hat W_*' (\pi^{(n)}, \lambda^{(n)}, \tau^{(n)})
(q_{\alpha_b}) |\ge  0\, .
$$

\vskip .3 truecm\noindent\proof We have
$$
\hat W_*'(\pi , \lambda , \tau )  = \hat V_*'(\pi ,\tau ) - {<\lambda ,v>\over <\lambda ,q>}\, .
$$
Therefore one has $|\hat W_*' (\pi , \lambda , \tau )(q_{\alpha_b}) |<\varepsilon$ if and only if
$$
\left| {<\lambda ,v>\over <\lambda ,q>}- \hat V_*'(\pi ,\tau ) (q_{\alpha_b})\right| <\varepsilon\, .
$$
For fixed $(\pi ,\tau )$, the set of $\lambda $ such that $|\hat W_*' (\pi , \lambda , \tau ) (q_{\alpha_b})|<\varepsilon $ has therefore Lebesgue measure  at most $C\varepsilon $ (because
$q$ and $v$ are normalized to have $l^2$ norm $1$, and $q$ is positive). Then, the Zorich invariant
measure of the same set is at most $C'\varepsilon^{1/(d-1)}$ (indeed from the estimate (6.4), p. 433 in [Y1], it follows that the Lebesgue measure of the set where the Zorich density is $\ge 2^m$ is at most
$C_02^{-m(1+{1\over d-2})}$, then, any set of Lebesgue measure $\varepsilon$ has Zorich measure at most $C_1\varepsilon^{1/(d-1)}$). The Claim now follows by a Borel--Cantelli argument. \qed

\vskip .5 truecm \noindent {\bf 3.7.5 End of the proof of Proposition 3.3.1 }
We have just seen that the quantity at the end of Section 3.7.3
$$
\Theta_2^{(n)} (W_{*}(y^\dagger)-W_{*}(x_{*}^{\hbox{\sevenrm max}}))
$$
(with $y^\dagger,x_{*}^{\hbox{\sevenrm max}}$ separated by $q_{\alpha_b}$ )  grows
in absolute value exponentially fast at rate $\theta_2$ in Zorich time $Z(n)$. This quantity
was seen in Section 3.7.3 to control $S_jw(x^*)$ for $Q^+(n)<j\le Q^+(n+1)$ (in the case $x_{*}^{\hbox{\sevenrm max}}<
q_{\alpha_b}$).

On the other hand, as $Q^+(n+1)$ is controlled by the return times in $I^{(n+1)}$,
we will have by the scaling rules
$$
\limsup_{n\rightarrow +\infty}{1\over Z(n)}\log Q^+(n+1) =\theta_1\; .
$$
(Actually, it is easy to see that almost surely $Q^+(n)$ grows exactly at rate $\theta_1$ in Zorich time).
The proof of Proposition 3.3.1 is now complete. \qed

%

\vskip .5 truecm\noindent{\bf 3.7.6. Remark.} Let $r=r(\pi,\lambda, \tau)$ be
the number of sequences $(\alpha_n)_{n\ge 0}$ such that $\eta_n(\alpha_{n+1})=
\alpha_n$ for all $n\ge 0$. It is clear that $r$ is invariant under the Rauzy--Veech
dynamics and therefore is constant almost everywhere. We claim that the dimension
of the simplex $\hbox{Aff}^{(1)}\, (\underline\gamma , w)$ is $r-1$; more precisely,
in the notation of Section 2.2, we claim that $r_d=r$ and $r_c=0$.

To prove this, we first observe that, to the $r$ sequences $(\alpha_n)_{n\ge 0}$, we
can associate points $x_1^*, x_2^*,\ldots ,x_r^*$ belonging to distinct orbits
and satisfying the estimate of Proposition 3.3.1. Then, we blow up these $r$ orbits as
in Section 3.3.2 to obtain an affine i.e.m.\ in $\hbox{Aff}^{(1)}\, (\underline\gamma , w)$
with $r$ orbits of wandering intervals. This shows that $r_d\ge r$.

Conversely, let $T_0\in\hbox{Aff}^{(1)}\, (\underline\gamma , w)$. We will show that there
are at most $r$ orbits of maximal wandering intervals and that the union of these
orbits has full measure. This implies $r_d\le r$ and $r_c=0$.

Let $h$ be the semiconjugacy from $T_0$ to $T_{\pi,\lambda}$. After changing names
if necessary, we can assume that there exists $0\le r'\le r$ such that
$h^{-1}(x_l^*)$ is a non trivial (wandering) interval for $1\le l\le r'$ and
a point for $r'<l\le r$. Denote by $(\alpha_n^{(l)})_{n\ge 0}$ the sequence
associated to $x_l^*$.

It follows from the estimates in Section 3.7 that there exists a constant
$C=C(\pi,\lambda,\tau)$ such that, for all $\alpha\in\cA$, $n\ge 0$
$$
\sum_{0\le i<Q_\alpha (n)}\exp(S_\alpha w(i))\le C
\hbox{Max}_{0\le i<Q_\alpha (n)}\exp (S_\alpha w(i))\, ,
$$
where $S_\alpha w(i)$ is as in 3.7.2 and $Q_\alpha (n)$ is the return time of
$I_\alpha^{(n)}$ in $I^{(n)}$.

From this we get
$$
\sum_{0\le i<Q_\alpha (n)}|T_0^i(h^{-1}(I_\alpha^{(n)}))|\le C
|h^{-1}(I_\alpha^{\hbox{\sevenrm max}}(n))|\, .
$$
Fix $\varepsilon >0$. Take then $n_0$ large enough to have the
$\alpha_{n_0}^{(l)}$ ($1\le l\le r$) distinct and
$$
|h^{-1}(I_{\alpha_{n_0}^{(l)}}^{\hbox{\sevenrm max}}(n_0))|\le
|h^{-1}(x_l^*)|+\varepsilon\, .
$$
Next, take $n_1>n_0$ large enough to have that the image by $\eta_{n_0}\circ\ldots
\circ\eta_{n_1-1}$ of $\cA$ is $\{\alpha_{n_0}^{(l)}\, , \, 1\le l\le r\}$.

The intervals $T_0^i(h^{-1}(I_\alpha^{(n_1)}))$, for $\alpha\in\cA$, $0\le i<Q_\alpha
(n_1)$, from a partition $\hbox{mod}\, 0$ of $[0,1]$.

If $\alpha$ is not one of the $\alpha_{n_1}^{(l)}$, writing
$\eta_{n_0}\circ\ldots\circ\eta_{n_1-1}(\alpha )=\alpha_{n_0}^{(j)}$, we have
that $h^{-1}(I_\alpha^{\hbox{\sevenrm max}}(n_1))\subset h^{-1}(I_{\alpha_{n_0}^{(j)}}^{\hbox{\sevenrm max}}
(n_0)\setminus I_{\alpha_{n_1}^{(j)}}^{\hbox{\sevenrm max}}
(n_1))$ hence $|h^{-1}(I_{\alpha}^{\hbox{\sevenrm max}}
(n_1))|<\varepsilon$.
When $\alpha=\alpha_{n_1}^{(l)}$ with $r'<l\le r$, we have also
$|h^{-1}(I_{\alpha}^{\hbox{\sevenrm max}}
(n_1))|<\varepsilon$.

In both cases, we obtain from above:
$$
\sum_{0\le i<Q_\alpha (n_1)}|T_0^i(h^{-1}(I_\alpha^{(n)}))|\le C\epsilon\, .
$$
Finally, when $\alpha$ is one of the $\alpha_{n_1}^{(l)}$, $1\le l\le r'$, we have
$$
\sum_{0\le i<Q_\alpha (n_1)}|T_0^i(h^{-1}(I_\alpha^{(n)})\setminus J_l^{(n)})|
\le C\varepsilon\, ,
$$
where $J_l^{(n)}$ is the wandering interval such that $T_0^m(J_l^{(n)})=h^{-1}(x_l^*)$
for some $0\le m<Q_\alpha (n_1)$. We have proved that the measure of the complement
of the orbits of the wandering intervals $h^{-1}(x_l^*)$ is arbitrarily small. This proves our claim.
\qed

\vfill\eject

\noindent {\bf References}

\vskip .3 truecm \noindent

\item{[AV]} A. Avila, M. Viana ``Simplicity of Lyapunov spectra: proof of the Zorich-Kontsevich conjecture''
{\it Acta Mathematica} {\bf 198} (2007) 1--56
\item{[AGY]} A. Avila, S. Gou\"ezel and J.-C. Yoccoz ``Exponential mixing for the Teichm\"uller flow''
preprint (2005)
\item{[BHM]} X. Bressaud, P. Hubert and A. Maass ``Persistence of wandering intervals in self-similar affine interval exchange transformations'' preprint arXiv: [math.DS] 0801.2088v1
\item{[CG]} R. Camelier, C. Gutierrez ``Affine interval exchange transformations with
wandering intervals'' {\it Ergod. Th. Dyn. Sys. } {\bf 17} (1997) 1315--1338
\item{[Co]} M. Cobo ``Piece-wise affine maps conjugate to
interval exchanges'' {\it Ergod. Th. Dyn. Sys. } {\bf 22} (2002) 375--407
\item{[Fo]} G. Forni ``Deviation of ergodic averages for area-preserving
flows on surfaces of higher genus.'' {\it Annals of Mathematics}
{\bf 155 } (2002) 1--103.
\item{[He]} M.~R. Herman ``Sur la conjugaison diff\'erentiable
des diff\'eomorphismes du cercle \`a des rotations'' {\it Publ. Math.
IHES} {\bf 49} (1979) 5--233
\item{[Ka]} A.B. Katok, ``Invariant measures of flows on orientable surfaces'' {\it Dokl. Akad. Nauk SSSR}
{\bf 211} (1973) 775--778
\item{[Ke]} M. Keane ``Interval exchange transformations''
{\it Math. Z.} {\bf 141} (1975) 25--31
\item{[L]} G. Levitt ``La d\'ecomposition dynamique et la diff\'erentiabilit\'e des feuilletages des surfaces''
{\it Ann. Inst. Fourier} {\bf 37} (1987) 85--116
\item{[LM]} I. Liousse, H. Marzougui ``\'Echanges d'intervalles affines conjugu\'es \`a des lin\'eaires''
{\it Ergod. Th. Dyn. Sys.} {\bf 22} (2002) 535--554
\item{[Ma]} H. Masur ``Interval exchange transformations and measured
foliations'' {\it Annals of Mathematics} {\bf 115} (1982) 169--200
\item{[MMY]} S. Marmi, P. Moussa and J.--C. Yoccoz ``The cohomological equation for Roth-type interval exchange maps''
{\it  J. Amer. Math. Soc.} {\bf  18} (2005)  823--872
\item{[MSS]} R. Mane, P. Sad and D. Sullivan ``On the dynamics of rational maps''
{\it Ann. Sci. E.N.S.} {\bf 16} (1983) 193--217
\item{[Ra]} G. Rauzy ``\'Echanges d'intervalles et transformations
induites'' {\it Acta Arit.} (1979) 315--328
\item{[V1]} W. Veech ``Interval exchange transformations''
{\it Journal d'Analyse Math\'e\-ma\-tique} {\bf 33} (1978) 222--272
\item{[V2]} W. Veech ``Gauss measures for transformations on the space of
interval exchange maps'' {\it Ann. of Math.} {\bf 115} (1982) 201--242
\item{[Y1]} J.--C. Yoccoz ``Continued fraction algorithms for interval exchange maps: an
introduction'', in Frontiers in Number Theory, Physics and Geometry I, Cartier P., Julia B., Moussa P.
and Vanhove P. editors, Springer--Verlag (2006)
\item{[Y2]} J.--C. Yoccoz ``Echanges d'intervalles'' Cours Coll\`ege de France 2005,   \par\noindent
\phantom{\item{[Y2]}} http://www.college-de-france.fr/media/equ$_{-}$dif/UPL8726$_{-}$yoccoz05.pdf
\item{[Z1]} A. Zorich ``Finite Gauss measure on the space of interval
exchange transformations. Lyapunov exponents''
Annales de l'Institut Fourier Tome
46 fasc. 2 (1996) 325-370
\item{[Z2]} A. Zorich ``Deviation for interval exchange
transformations'' {\it Ergod. Th. Dyn. Sys.}{\bf 17} (1997),
1477--1499

\bye